\documentclass[a4paper,12pt]{amsart}
\usepackage{amssymb, amsmath, graphicx}
\usepackage[curve]{xypic}

\providecommand{\id}{\textnormal{id}}

\providecommand{\Ker}{\textnormal{Ker}}
\providecommand{\IIm}{\textnormal{Im}}
\providecommand{\Coker}{\textnormal{Coker}}

\providecommand{\Hom}{\textnormal{Hom}}

\providecommand{\Ker}{\textnormal{Ker}}

\providecommand{\ch}{\textnormal{ch}}

\providecommand{\cpt}{\textnormal{cpt}}

\providecommand{\fl}{\textnormal{fl}}

\providecommand{\N}{\mathbb{N}}
\providecommand{\Z}{\mathbb{Z}}

\providecommand{\R}{\mathbb{R}}
\providecommand{\C}{\mathbb{C}}
\providecommand{\h}{\mathfrak{h}}
\providecommand{\kk}{\mathfrak{k}}
\providecommand{\cl}{\textnormal{cl}}
\providecommand{\op}{\textnormal{op}}
\providecommand{\dR}{\textnormal{dR}}
\providecommand{\Td}{\textnormal{Td}}

\allowdisplaybreaks

\begin{document}

\title{Flat pairing and generalized Cheeger-Simons characters}
\author{Fabio Ferrari Ruffino}
\address{Departamento de Matem\'atica - Universidade Federal de S\~ao Carlos - Rod.\ Washington Lu\'is, Km 235 - C.P.\ 676 - 13565-905 S\~ao Carlos, SP, Brasil}
\email{ferrariruffino@gmail.com}
\thanks{I would like to thank the referee for many helpful suggestions, critics and remarks, that changed the structure of the paper. I was supported by FAPESP (Funda\c{c}\~ao de Amparo \`a Pesquisa do Estado de S\~ao Paulo), processo 2014/03721-3, and I am thankful to the ICMC - USP, S\~ao Carlos, for its support during my staying there.}

\begin{abstract}
Let $h^{\bullet}$ be a multiplicative cohomology theory, $h_{\bullet}$ its dual homology theory and $\hat{h}^{\bullet}$ a differential refinement. We first construct the natural pairing between $h_{\bullet}$ and the flat part of $\hat{h}^{\bullet}$, generalizing the holonomy of a flat Deligne cohomology class. Then, in order to generalize the holonomy of any Deligne cohomology class, we define the generalized Cheeger-Simons characters. The latter are functions from suitably defined differential cycles to the cohomology ring of the point, such that the value on a trivial cycle only depends on the curvature.
\end{abstract}

\maketitle

\newtheorem{Theorem}{Theorem}[section]
\newtheorem{Lemma}[Theorem]{Lemma}
\newtheorem{Corollary}[Theorem]{Corollary}
\newtheorem{ThmDef}[Theorem]{Theorem - Definition}

\theoremstyle{definition}
\newtheorem{Rmk}[Theorem]{Remark}
\newtheorem{Rmks}[Theorem]{Remarks}
\newtheorem{Def}[Theorem]{Definition}

\section{Introduction}

Let us consider the ordinary differential cohomology $\hat{H}^{\bullet}$ on a smooth manifold $X$ \cite{HS}. The group $\hat{H}^{n}(X)$ is canonically isomorphic to the group of the Cheeger-Simons differential characters of degree $n$. An element of the latter is a couple $(\chi, \omega)$, where $\chi$ is an $\R/\Z$-valued group morphism defined on the smooth $(n-1)$-cycles of $X$ (whose exponential is the holonomy), and $\omega$ is an integral $n$-form on $X$ (the curvature) such that, on a $p$-boundary $\partial D$, one has:
	\[\chi(\partial D) = \int_{D} \omega \mod \Z.
\]
In particular, when the class is flat, the holonomy only depends on the homology class of the cycle, and actually the flat part of $\hat{H}^{n}(X)$ is canonically isomorphic to $H^{n-1}(X; \R/\Z)$.

The aim of the present paper is to generalize this picture to a differential refinement $\hat{h}^{\bullet}$ of any multiplicative cohomology theory $h^{\bullet}$ \cite{Bunke, HS, Upmeier}. We start considering the flat case: a differential class of degree $n$ provides a morphism from the homology theory $h_{\bullet}$ to the cohomology ring of the point; we thus get a pairing between $h_{\bullet}$ and the flat part of $\hat{h}^{\bullet}$. Then, in order to consider even non-flat classes, we consider the geometrical definition of the homology theory $h_{\bullet}$ dual to $h^{\bullet}$, as described in \cite{Jakob}, and we define a variant of that construction, using a suitable differential extension of cycles and boundaries. In this way it is possible to define the the generalized Cheeger-Simons characters as functions from the differential cycles to the cohomology ring of the point, such that the value on a trivial cycle only depends on the curvature. In the case of $K$-theory this construction is equivalent to the one defined in \cite{BM}.

The paper is organized as follows. In sections \ref{DiffCoh} and \ref{OrInt} we recall the preliminaries about differential cohomology and the corresponding notions of orientation and integration. In section \ref{CohHom} we recall the notion of dual homology theory and its basic properties. In section \ref{FlatPairing} we show the pairing between $h_{\bullet}$ and the flat part of $\hat{h}^{\bullet}$. Finally, in section \ref{DiffCycles} we construct a model of the homology groups via differential cycles and we define the generalized Cheeger-Simons characters.

\section{Differential cohomology}\label{DiffCoh}

Let $\mathcal{M}$ be the category of smooth manifolds or of smooth compact manifolds (even with boundary), and let $\mathcal{A}_{\Z}$ be the category of $\Z$-graded abelian groups. We consider a cohomology theory $h^{\bullet}$, defined on a category including $\mathcal{M}$. We use the following notation:
	\[\h^{\bullet} := h^{\bullet}(\{pt\}) \qquad \h^{\bullet}_{\R} := \h^{\bullet} \otimes_{\Z} \R.
\]
Moreover, for any object $M$ of $\mathcal{M}$, we call $\ch: h^{\bullet}(M) \rightarrow H^{\bullet}_{\dR}(M; \h^{\bullet}_{\R})$ the generalized Chern character \cite[sec.\ 4.8 p.\ 47]{HS}. In this introductory section, we follow \cite[sec.\ 1]{BS}.
\begin{Def} A \emph{differential extension} of $h^{\bullet}$ is a functor $\hat{h}^{\bullet}: \mathcal{M}^{\op} \rightarrow \mathcal{A}_{\Z}$, together with the following natural transformations of $\mathcal{A}_{\Z}$-valued functors:
\begin{itemize}
	\item $I: \hat{h}^{\bullet}(M) \rightarrow h^{\bullet}(M)$;
	\item $R: \hat{h}^{\bullet}(M) \rightarrow \Omega_{\cl}^{\bullet}(M; \h^{\bullet}_{\R})$, called \emph{curvature};
	\item $a: \Omega^{\bullet-1}(M; \h^{\bullet}_{\R})/\IIm(d) \rightarrow \hat{h}^{\bullet}(M)$,
\end{itemize}
such that:
\begin{itemize}
	\item $R \circ a = d$;
	\item the following diagram is commutative:
		\[\xymatrix{
		\hat{h}^{\bullet}(M) \ar[r]^{I} \ar[d]_{R} & h^{\bullet}(M) \ar[d]^{\ch} \\
		\Omega_{\cl}^{\bullet}(M; \h^{\bullet}_{\R}) \ar[r]^(.47){dR} & H^{\bullet}_{\dR}(M; \h^{\bullet}_{\R});
	}\]
	\item the following sequence is exact:
		\[\xymatrix{
		\hat{h}^{\bullet-1}(M) \ar[r]^(.35){\ch} & \Omega^{\bullet-1}(M; \h^{\bullet}_{\R})/\IIm(d) \ar[r]^(.65){a} & \hat{h}^{\bullet}(M) \ar[r]^{I} & h^{\bullet}(M) \ar[r] & 0.
	}\]
\end{itemize}
We also call $\hat{h}^{\bullet}$ \emph{differential cohomology theory}.
\end{Def}
A class $\hat{\alpha} \in \hat{h}^{n}(M)$ is called \emph{flat} if and only if $R(\hat{\alpha}) = 0$. Considering flat classes, we get the functor $\hat{h}^{\bullet}_{\fl}: \mathcal{M} \rightarrow \mathcal{A}_{\Z}$. Thus, we get the following commutative hexagon:
\begin{equation}\label{CommHex1}
	\resizebox{0.9\textwidth}{!}{
	\xymatrix{
	& \Omega^{\bullet-1}(M; \h^{\bullet}_{\R})/\IIm(d) \ar[rr]^{d} \ar[dr]^{a} & & \Omega_{\cl}^{\bullet}(M; \h^{\bullet}_{\R}) \ar[dr]^{dR} \\
	H^{\bullet-1}_{\dR}(M; \h^{\bullet}_{\R}) \ar[ur] \ar[dr]^{a} & & \hat{h}^{\bullet}(M) \ar[ur]^{R} \ar[dr]^{I} & & H^{\bullet}_{\dR}(M; \h^{\bullet}_{\R}). \\
	& \hat{h}^{\bullet}_{\fl}(M) \ar@{^(->}[ur] \ar[rr]^{I} & & h^{\bullet}(M) \ar[ur]^{\ch}
}}
\end{equation}
The following lemma easily follows from the previous axioms.
\begin{Lemma}[Homotopy formula] If $\hat{\alpha} \in \hat{h}^{\bullet}(I \times X)$ and $i_{0}, i_{1}: X \rightarrow I \times X$ are the natural embeddings, we have:
\begin{equation}\label{HomFormula}
	i_{1}^{*}\hat{\alpha} - i_{0}^{*}\hat{\alpha} = a\biggl( \int_{I} R(\hat{\alpha}) \biggr).
\end{equation}
\end{Lemma}
We now introduce integration. Given a functor $\mathcal{F}: \mathcal{M} \rightarrow \mathcal{C}$, for any category $\mathcal{C}$, we define the functor $S\mathcal{F}: \mathcal{M} \rightarrow \mathcal{C}$ defined by $S\mathcal{F}(M) := \mathcal{F}(S^{1} \times M)$ on objects and $S\mathcal{F}(f) := \mathcal{F}(\id_{S^{1}} \times f)$ on morphisms.
\begin{Def} A \emph{differential extension with integration} of $h^{\bullet}$ is a differential extension $(\hat{h}^{\bullet}, I, R, a)$ together with a natural transformation:
	\[\int_{S^{1}}: S\hat{h}^{\bullet+1} \rightarrow \hat{h}^{\bullet},
\]
such that:
\begin{itemize}
	\item $\int_{S^{1}} \circ (t \times \id)^{*} = -\int_{S^{1}}$, where $t: S^{1} \rightarrow S^{1}$ is defined by $t(e^{i\theta}) := e^{-i\theta}$;
	\item $\int_{S^{1}} \circ p^{*} = 0$, where $p: S^{1} \times X \rightarrow X$ is the projection;
	\item the following diagram is commutative:
		\[\resizebox{0.8\textwidth}{!}{
\xymatrix{
		S\Omega^{\bullet}(M; \h^{\bullet}_{\R})/\IIm(d) \ar[r]^(.6){Sa} \ar[d]^{\int_{S^{1}}} & S\hat{h}^{\bullet+1}(M) \ar[r]^{SI} \ar[d]^{\int_{S^{1}}} \ar@/^2pc/[rr]^{SR} & Sh^{\bullet+1}(M) \ar[d]^{\int_{S^{1}}} & S\Omega_{\cl}^{\bullet+1}(M; \h^{\bullet}_{\R}) \ar[d]^{\int_{S^{1}}} \\
		\Omega^{\bullet-1}(M; \h^{\bullet}_{\R})/\IIm(d) \ar[r]^(.65){a} & \hat{h}^{\bullet}(M) \ar[r]^{I} \ar@/_2pc/[rr]_{R} & h^{\bullet}(M) & \Omega_{\cl}^{\bullet}(M; \h^{\bullet}_{\R}).
	}}\]
\end{itemize}
\end{Def}
Let us consider a differential extension with integration $(\hat{h}^{\bullet}, I, R, a, \int_{S^{1}})$. It is shown in \cite[pp.\ 27-32]{BS} that, if $h^{\bullet}$ is rationally even (i.e., $\h^{2k+1}_{\R} = 0$ for every $k \in \Z$), $\h^{k}$ is countably generated for every $k \in \Z$ and $\mathcal{M}$ is the category of all smooth manifolds, there is an isomorphism of functors $\hat{h}^{\bullet}_{\fl} \simeq h^{\bullet-1}(\,\cdot\,; \R/\Z)$. When $\mathcal{M}$ is the category of compact manifolds, we must require that $\h^{2k+1} = 0$ and $\h^{2k}$ is finitely generated for every $k \in \Z$. In these cases, the commutative hexagon \eqref{CommHex1} becomes:
\begin{equation}\label{CommHex2}
	\resizebox{0.8\textwidth}{!}{
	\xymatrix{
	& \Omega^{\bullet-1}(M; \h^{\bullet}_{\R})/\IIm(d) \ar[rr]^{d} \ar[dr]^{a} & & \Omega_{\cl}^{\bullet}(M; \h^{\bullet}_{\R}) \ar[dr]^{dR} \\
	H^{\bullet-1}_{\dR}(M; \h^{\bullet}_{\R}) \ar[ur] \ar[dr] & & \hat{h}^{\bullet}(M) \ar[ur]^{R} \ar[dr]^{I} & & H^{\bullet}_{\dR}(M; \h^{\bullet}_{\R}), \\
	& h^{\bullet-1}(M; \R/\Z) \ar@{^(->}[ur] \ar[rr]^(.55){\beta} & & h^{\bullet}(M) \ar[ur]^{\ch}
}}
\end{equation}
where $\beta$, in the last line, is the Bockstein map of the long exact sequence induced by the coefficient sequence $0 \rightarrow \Z \rightarrow \R \rightarrow \R/\Z \rightarrow 0$.

Finally, we introduce products, thus we suppose that $h^{\bullet}$ is a \emph{multiplicative} cohomology theory. We call $\mathcal{R}_{\Z}$ the category of $\Z$-graded commutative rings. There is a natural forgetful functor $\mathcal{R}_{\Z} \rightarrow \mathcal{A}_{\Z}$, that we apply when needed, without writing it explicitly.
\begin{Def}\label{MultDiffExt} A \emph{multiplicative differential extension} of $h^{\bullet}$ is a differential extension $(\hat{h}^{\bullet}, I, R, a)$ such that $\hat{h}^{\bullet}: \mathcal{M}^{\op} \rightarrow \mathcal{R}_{\Z}$ and:
\begin{itemize}
	\item $I$ and $R$ are multiplicative;
	\item $\hat{\alpha} \cdot a(\omega) = a(R(\hat{\alpha}) \wedge \omega)$ for every $\hat{\alpha} \in \hat{h}^{\bullet}(M)$ and $\omega \in \Omega^{\bullet}(M; \h^{\bullet}_{\R})/\IIm(d)$.
\end{itemize}
\end{Def}
In the following we consider multiplicative differential extensions with integration.

\section{Orientation and integration}\label{OrInt}

Following \cite[sec.\ 4.8-4.10]{Bunke}, we briefly recall the topological notions of orientation and integration and we extend them to the differential case.

\subsection{Topological orientation and integration}

Let $X$ be a compact manifold and $E \rightarrow X$ a real vector bundle of rank $n$. We call $E_{0}$ the fiber sub-bundle of $E$, containing the complement of the unit ball of each fiber of $E$, with respect to any fixed metric. The bundle $E$ is orientable with respect to a multiplicative cohomology theory $h^{\bullet}$ if there exists a \emph{Thom class} $u \in h^{n}(E, E_{0})$ \cite[p.\ 253]{Rudyak}, i.e.\ a class such that, for every $x \in X$, the restriction $u\vert_{(E_{x}, E_{0,x})}$ maps to $\pm 1$ under the isomorphism:
	\[h^{n}(E_{x}, E_{0,x}) \simeq h^{n}(\R^{n}, \R^{n} \setminus B^{n}) \simeq \tilde{h}^{n}(S^{n}) \simeq \mathfrak{h}^{0}.
\]
Using the product $h^{\bullet}(E, E_{0}) \otimes_{\Z} h^{\bullet}(E) \rightarrow h^{\bullet}(E, E_{0})$, we define the Thom isomorphism $\alpha \mapsto u \cdot \pi^{*}\alpha$, between $h^{\bullet}(X)$ and $h^{\bullet+n}(E, E_{0})$. Thus, we define the integration map $\int_{E/X} \colon h^{\bullet}(E, E_{0}) \rightarrow h^{\bullet-n}(X)$ by $u \cdot \pi^{*}\alpha \mapsto \alpha$. If the characteristic of $\h^{\bullet}$ is different from $2$, the $n$-degree component of $\ch \,u$ defines an orientation of $E$ in the usual sense, hence it is possible to integrate a compactly-supported form fiber-wise. We define the \emph{Todd class} $\Td(u) := \int_{E/X} \ch \, u \in H^{0}_{\dR}(X; \h^{\bullet}_{\R})$. The following formula holds:
\begin{equation}\label{ChIntegral}
	\int_{E/X} \ch \, \alpha = \Td(u) \cdot \biggl( \ch \int_{E/X} \alpha \biggr).
\end{equation}
\begin{Lemma}[2x3 principle]\label{Rule23Top} Given two bundles $E, F \rightarrow X$, with projections $p_{E}: E \oplus F \rightarrow E$ and $p_{F}: E \oplus F \rightarrow F$, we consider a triple $(u, v, w)$ of Thom classes on $E$, $F$ and $E \oplus F$ respectively, such that $w = p_{E}^{*}u \cdot p_{F}^{*}v$. Two elements of such a triple uniquely determines the third one.
\end{Lemma}
For the proof see \cite[prop.\ 1.10 p.\ 307]{Rudyak}. In order to define the Gysin map associated to a smooth neat\footnote{We recall that $f: Y \rightarrow X$ is \emph{neat} if $f^{-1}(\partial X) = \partial Y$ and $df_{y}: T_{y}Y/T_{y}(\partial Y) \rightarrow T_{f(x)}X/T_{f(x)}(\partial X)$ is an isomorphism for every $y \in \partial Y$ \cite[Appendix C.2 p.\ 90]{HS}. See also \cite[pp.\ 27ff]{Kosinski} about neat submanifolds and tubular neighborhoods.} map of compact manifolds $f: Y \rightarrow X$, we consider a neat proper embedding $\iota: Y \hookrightarrow X \times \R^{N}$, such that $\pi_{X} \circ \iota = f$, and we endow the normal bundle $N_{\iota(Y)}(X \times \R^{N})$ with a Thom class $u$. Then, we consider a neat tubular neighborhood $U$ of $\iota(Y)$ in $X \times \R^{N}$, a diffeomorphism $\varphi: N_{\iota(Y)}(X \times \R^{N}) \rightarrow U$ and the natural inclusion $i: U \rightarrow X$, inducing a push-forward in compactly-supported cohomology. The Gysin map $f_{!}: h^{\bullet}(Y) \rightarrow h^{\bullet - n}(X)$, being $n = \dim Y - \dim X$, is defined as:
\begin{equation}\label{GysinEmbeddingTop}
	f_{!}(\alpha) = \int_{\R^{N}}i_{*}\varphi_{*}(u \cdot \pi^{*}\alpha).
\end{equation}
We used the fact that, since $X$ is compact, we have $h^{\bullet}(E, E_{0}) \simeq h^{\bullet}_{\cpt}(E)$. The integration map is defined as follows: since $h^{\bullet}_{\cpt}(X \times \mathbb{R}^{N}) = \tilde{h}^{\bullet}((X \times \mathbb{R}^{N})^{+}) \simeq \tilde{h}^{\bullet}(\Sigma^{N}(X_{+}))$, for $X_{+} = X \sqcup \{\infty\}$, we apply the suspension isomorphism $\tilde{h}^{\bullet}(\Sigma^{N}(X_{+})) \simeq \tilde{h}^{\bullet-n}(X_{+}) \simeq h^{\bullet-n}(X)$.
This construction of the Gysin map naturally leads to the following definition.
\begin{Def}\label{DiffOrientedMap} A \emph{representative of an $h^{\bullet}$-orientation} of a smooth neat map between compact manifolds $f: Y \rightarrow X$ is the datum of:
\begin{itemize}
	\item a neat embedding $\iota: Y \hookrightarrow X \times \R^{N}$, for any $N \in \N$, such that $\pi_{X} \circ \iota = f$;
	\item a Thom class $u$ of the normal bundle $N_{\iota(Y)}(X \times \R^{N})$;
	\item a diffeomorphism $\varphi: N_{\iota(Y)}(X \times \R^{N}) \rightarrow U$, for $U$ a neat tubular neighborhood of $\iota(Y)$ in $X \times \R^{N}$.
\end{itemize}
\end{Def}
Let us consider a vector $v_{y} \in N_{\iota(Y)}(X \times \R^{N})_{\iota(y)}$. It is sent by $\varphi$, as defined in \ref{DiffOrientedMap}, to a point $\varphi(v_{y}) \in X \times \R^{N}$. When $f$ is a submersion, we can require that the first component of $\varphi(v_{y})$ is $f(y)$. This means that the following diagram commutes:
\begin{equation}\label{PropRepres}
	\xymatrix{
	N_{\iota(Y)}(X \times \R^{N}) \ar[r]^(.7){\varphi} \ar[d]_{\pi_{N}} & U \ar[d]^{\pi_{X}} \\
	\iota(Y) \ar[r]^{\pi_{X}} & X.
}\end{equation}
\begin{Def}\label{PropDef} A representative of an $h^{\bullet}$-orientation of $f: Y \rightarrow X$ is \emph{proper} if diagram \eqref{PropRepres} commutes.
\end{Def}
\begin{Lemma}\label{SubProp} Let $f: Y \rightarrow X$ be a submersion. Then, for any neat embedding $\iota: Y \hookrightarrow X \times \R^{N}$ and any Thom class $u$ of the normal bundle, there exists a \emph{proper} representative $(\iota, u, \varphi)$ of an $h^{\bullet}$-orientation of $f$.
\end{Lemma}
\begin{proof} We define the following projection:
	\[\begin{split}
	p: \; &T(X \times \R^{N})\vert_{\iota(Y)} \rightarrow (\iota(Y) \times \R^{N})/d\iota(\Ker \, df) \\
	& (df_{y}(v), w) \mapsto [w - \pi_{\R^{N}}(d\iota_{y}(v))].
\end{split}\]
It is easy to prove that $p$ is surjective and $\Ker \, p = T(\iota(Y))$, hence $N_{\iota(Y)}(X \times \R^{N}) \simeq (\iota(Y) \times \R^{N})/d\iota(\Ker \, df)$. Using the standard metric on $\R^{N}$, we can identify the quotient $(\iota(Y) \times \R^{N})/d\iota(\Ker \, df)$ with the sub-bundle $d\iota(\Ker \, df_{y})^{\bot}$, hence we can apply $\varphi$ in such a way that the fiber on $\iota(y)$ is sent to an open subset of $\{f(y)\} \times \R^{N}$ (the reader can work out the details). \end{proof}

We now introduce a suitable equivalence relation among representatives of orientations. First we need a generalization of definition \ref{PropDef}. Let us consider a representative $(J, U, \Phi)$ of an $h^{\bullet}$-orientation of $\id \times f: I \times Y \rightarrow I \times X$ and a neighborhood $V \subset I$ of $\{0, 1\}$. We say that the representative is \emph{proper on $V$} if a vector $(x,v)_{(t,y)} \in N_{\iota(V \times Y)}(V \times X \times \R^{N})_{\iota(t,y)}$ is sent by $\Phi$ to a point $\Phi((x,v)_{(t,y)}) \in V \times X \times \R^{N}$ whose first component is $t$. This means that the following diagram commutes:
\begin{equation}\label{PropRepres2}
	\xymatrix{
	N_{\iota(V \times Y)}(V \times X \times \R^{N}) \ar[r]^(.77){\Phi} \ar[d]_{\pi_{N}} & U \ar[d]^{\pi_{I}} \\
	\iota(V \times Y) \ar[r]^(.6){\pi_{I}} & I.
}\end{equation}
In this case, calling $f_{0} := \id_{\{0\}} \times f$ and $f_{1} := \id_{\{1\}} \times f$, we can define the restrictions $(J, U, \Phi)\vert_{f_{0}}$ and $(J, U, \Phi)\vert_{f_{1}}$.

\begin{Def}\label{HomotopyOrientations} A \emph{homotopy} between two representatives $(\iota, u, \varphi)$ and $(\iota', u', \varphi')$ of an $h^{\bullet}$-orientation of $f: Y \rightarrow X$ is a representative $(J, U, \Phi)$ of an $h^{\bullet}$-orientation of $\id \times f: I \times Y \rightarrow I \times X$, such that:
\begin{itemize}
	\item $(J, U, \Phi)$ is proper over a neighborhood $V \subset I$ of $\{0, 1\}$;
	\item $(J, U, \Phi)\vert_{f_{0}} = (\iota, u, \varphi)$ e $(J, U, \Phi)\vert_{f_{1}} = (\iota', u', \varphi')$.
\end{itemize}
\end{Def}
On the trivial bundle $X \times \R^{N}$ there is a canonical Thom class, defined in the following way. On $pt \times \R^{N}$, whose compactification is $pt \times S^{N}$, we put the class $u_{0} \in \tilde{h}^{N}(S^{N})$ corresponding to the suspension of $1 \in \h^{0}$. Then, we put on $X \times \R^{N}$ the class $\pi_{\R^{N}}^{*}u_{0}$.
\begin{Def}\label{EquivStab} Let us consider a representative $(\iota, u, \varphi)$, with $\iota: Y \hookrightarrow X \times \R^{N}$.
\begin{itemize}
	\item For any $L \in \N$, we define $\iota': Y \hookrightarrow X \times \R^{N + L}$ by $\iota'(y) := (\iota(y), 0)$. Then $N_{\iota'(Y)}(X \times \R^{N+L}) \simeq N_{\iota(Y)}(X \times \R^{N}) \oplus (\iota(Y) \times \R^{L})$.
	\item We put the canonical orientation $u_{0}$ on the trivial bundle $\iota(Y) \times \R^{L}$, and the orientation $u'$ induced by $u$ and $u_{0}$ on $N_{\iota'(Y)}(X \times \R^{N+L})$.
	\item Finally, for $v_{y} \in N_{\iota(Y)}(X \times \R^{N})$ and $w \in \R^{L}$, we define $\varphi'(v_{y}, w) := (\varphi(v_{y}), w) \in X \times \R^{N+L}$.
\end{itemize}
The representative $(\iota', u', \varphi')$ is called \emph{equivalent by stabilization} to $(\iota, u, \varphi)$. \end{Def}

\begin{Def} A \emph{$\hat{h}^{\bullet}$-orientation} on $f: Y \rightarrow X$ is an equivalence class $[\iota, u, \varphi]$ of representatives, up to the equivalence relation generated by homotopy and stabilization.
\end{Def}
The Gysin map $f_{!}$ only depends on the $\hat{h}^{\bullet}$-orientation of $f$, not on the specific representative (\cite[theorem 5.24 p.\ 233]{Karoubi}, \cite[sec.\ 4.9]{Bunke}). Moreover, because of the uniqueness up to homotopy and stabilization of the tubular neighborhood and of the diffeomorphism with the normal bundle, the class $[\iota, u, \varphi]$ does not depend on $\varphi$, hence we denote it by $[\iota, u]$.

\begin{Lemma}\label{GysinMapProp1} Given an $h^{\bullet}$-oriented map $f: Y \rightarrow X$, the Gysin map is a morphism of $\hat{h}^{\bullet}(X)$-modules, i.e., for any $\alpha \in h^{\bullet}(Y)$ and $\beta \in h^{\bullet}(X)$:
	\begin{equation}\label{GysinPullBack}
	f_{!}(\alpha \cdot f^{*}\beta) = f_{!}(\alpha) \cdot \beta.
\end{equation}
\end{Lemma}
For the proof see \cite[theorem 5.24 p.\ 233]{Karoubi}.
\begin{Def}\label{OrientedMapComposition} Let $f: Y \rightarrow X$ and $g: X \rightarrow W$ be $h^{\bullet}$-oriented maps, with orientations $[\iota, u]$ and $[\kappa, v]$, where $\iota: Y \hookrightarrow X \times \R^{N}$ and $\kappa: X \hookrightarrow W \times \R^{L}$. There is a naturally induced $h^{\bullet}$-orientation on $g \circ f: Y \rightarrow W$, that we denote by $[\kappa, v][\iota, u]$, defined in the following way:
\begin{itemize}
	\item we choose the embedding $\xi = (\kappa, \id_{\R^{N}}) \circ \iota: Y \hookrightarrow W \times \R^{L+N}$;
	\item on the normal bundle $N_{\xi(Y)}(W \times \R^{L+N}) \simeq N_{\iota(Y)}(X \times \R^{L}) \oplus N_{\kappa(X) \times \R^{L}}(W \times \R^{L+N})\vert_{\xi(Y)} \simeq N_{\iota(Y)}(X \times \R^{L}) \oplus (\pi^{*}_{N}N_{\kappa(X)}W \times \R^{L})\vert_{\xi(Y)}$, for $\pi_{L}: \R^{L+N} \rightarrow \R^{N}$, we put the Thom class $w$ induced from the ones on $N_{\iota(Y)}(X \times \R^{L})$ and $N_{\kappa(X)}(W \times \R^{N})$.
\end{itemize}
We set $[\kappa, v][\iota, u] := [\xi, w]$.
\end{Def}
\begin{Lemma}\label{GysinMapProp2} With the data of definition \ref{OrientedMapComposition}, $(g \circ f)_{!} = g_{!} \circ f_{!}$.
\end{Lemma}
For the proof see \cite[theorem 5.24 p.\ 233]{Karoubi}. The following lemma is a consequence of lemma \ref{Rule23Top} and of the uniqueness up to homotopy and stabilization of the embedding $\iota$.
\begin{Lemma}[2x3 principle]\label{Rule23TopMaps} Let $f: Y \rightarrow X$ and $g: X \rightarrow W$ be $h^{\bullet}$-oriented neat maps, with orientations $[\iota, u]$ and $[\kappa, v]$, and let $[\xi, w] := [\kappa, v][\iota, u]$ be the orientation induced on $g \circ f$. Two elements of the triple $([\iota, u], [\kappa, v], [\xi, w])$ uniquely determines the third one.
\end{Lemma}
Finally, we consider the orientation of manifolds.
\begin{Def}\label{OrientedManifold} An $h^{\bullet}$-orientation of a manifold without boundary $X$ is an $h^{\bullet}$-orientation of the map $p_{X}: X \rightarrow \{pt\}$.
\end{Def}
By definition, giving an orientation to $p_{X}$ means fixing an orientation $u$ on the (stable) normal bundle of $X$. We set $\Td(X) := \Td(u)$. If $Y$ and $X$ are oriented, because of the 2x3 principle a map $f: Y \rightarrow X$ inherits an orientation, hence the Gysin map is well-defined.
\begin{Def}\label{OrientedManifoldBoundary} An \emph{$h^{\bullet}$-orientation} on a smooth manifold \emph{with} boundary $X$ is the datum of:
\begin{itemize}
	\item a neat map $\Phi: X \rightarrow I$ such that $\partial X = \Phi^{-1}\{0\}$;
	\item an $h^{\bullet}$-orientation of $\Phi$.
\end{itemize}
\end{Def}
Again, giving an orientation to $\Phi$ means fixing an orientation $u$ on the stable normal bundle of $X$: in fact, since $\Phi^{-1}\{1\} = \emptyset$, the normal bundle in $I \times \R^{N}$ is canonically isomorphic to the the normal bundle in $[0, 1) \times \R^{N} \simeq \R^{N+1}_{+}$, being $\R^{N}_{+} = \{(x_{1}, \ldots, x_{N}) \in \R^{N} \,\vert\, x_{N} \geq 0\}$. We set $\Td(X) := \Td(u)$.

An orientation $[\iota, u]$ on a manifold with boundary canonically induces an orientation $[\iota', u']$ on the boundary: in fact, calling $i_{\partial X}: \partial X \hookrightarrow X$ the natural embedding, we set $\iota' := \iota \circ i_{\partial X}: \partial X \hookrightarrow \R^{N-1}$ and $u' := u\vert_{\partial X}$. The induced orientation does not depend on the representative chosen, since the tubular neighborhoods, the diffeomorphisms and the homotopies between them are neat, hence a homotopy between two representatives restricts to the boundary. The same holds for stabilization. Moreover, defining $\Phi: X \rightarrow I$ as in \ref{OrientedManifoldBoundary}, one has, for every $\alpha \in h^{\bullet}(X)$:
\begin{equation}\label{RestrictionBoundaryGysin}
	(p_{\partial X})_{!}(\alpha\vert_{\partial X}) = (\Phi_{!}\alpha)\vert_{\{0\}}.
\end{equation}
Such a formula is due to the fact that all the structures involved in the definition of the Gysin map for $p_{\partial X}$ are the restrictions to the boundary of the corresponding structures for $\Phi_{!}$.

\subsection{Differential orientation of a vector bundle}

If we consider a differential refinement $\hat{h}^{\bullet}$ of $h^{\bullet}$, in order to orient a vector bundle one just has to refine a Thom class $u$ to a \emph{differential Thom class}.
\begin{Def} Let $\hat{h}^{\bullet}$ be a multiplicative differential extension of $h^{\bullet}$. A \emph{differential Thom class} of $E$ is a compactly supported class $\hat{u} \in \hat{h}^{n}_{\cpt}(E)$ such that $I(\hat{u}) \in h^{n}_{\cpt}(E)$ is a Thom class for $h^{\bullet}$.
\end{Def}
Using the product $\hat{h}^{\bullet}_{\cpt}(E) \otimes_{\Z} \hat{h}^{\bullet}(E) \rightarrow \hat{h}^{\bullet}_{\cpt}(E)$, we define the differential Thom morphism, which is not surjective any more, as $\hat{\alpha} \mapsto \hat{u} \cdot \pi^{*}\hat{\alpha}$. We define the \emph{Todd class} $\Td(\hat{u}) := \int_{E/X} R(\hat{u}) \in \Omega^{0}_{\cl}(X; \h^{\bullet}_{\R})$. If follows that $I(\Td(\hat{u})) = \Td(I(\hat{u}))$.

\begin{Def}\label{HomotopyThom} Let $\pi_{X}: I \times X \rightarrow X$ be the natural projection and $i_{0}, i_{1}: X \rightarrow I \times X$ the natural embeddings. Two differential Thom classes $\hat{u}, \hat{u}' \in \hat{h}^{n}_{\cpt}(E)$ are \emph{homotopic} if there exists a Thom class $\hat{U} \in \hat{h}^{n}_{\cpt}(\pi_{X}^{*}E)$ such that $i_{0}^{*}\hat{U} = \hat{u}$, $i_{1}^{*}\hat{U} = \hat{u}'$ and $\Td(\hat{U}) = \pi_{X}^{*}\Td(\hat{u})$.
\end{Def}
Let us consider two Thom classes $\hat{u}, \hat{u}' \in \hat{h}^{n}_{\cpt}(E)$, refining $u \in h^{n}_{\cpt}(E)$. We have that $\hat{u}' - \hat{u} = a(\eta)$, with $\eta \in \Omega^{n-1}_{\cpt}(E; \mathfrak{h}^{\bullet}_{\R})/\IIm d$. Hence:
\begin{equation}\label{dIntEta}
	d\int_{E/X} \eta = \int_{E/X} d\eta = \int_{E/X} R(a(\eta)) = \int_{E/X} \bigl( R(\hat{u}') - R(\hat{u}) \bigr) = \Td(\hat{u}') - \Td(\hat{u}).
\end{equation}
If $\Td(\hat{u}') = \Td(\hat{u})$, then $d\int_{E/X} \eta = 0$ (this happens in particular if $\hat{u}$ and $\hat{u}'$ are homotopic). In this case, there is a well-defined de-Rham cohomology class:
\begin{equation}\label{dRIntEta}
	\dR\biggl(\int_{E/X} \eta\biggr) \in H^{-1}_{\dR}(X; \mathfrak{h}^{\bullet}_{\R}).
\end{equation}
Moreover, the form $\eta$ is unique up to a form $\xi$ representing a class belonging to the image of the Chern character. Because of formula \eqref{ChIntegral}, we have that $\dR(\xi) \in \IIm(\ch)$ if and only if $\dR(\int_{E/X} \xi) \in \Td(u) \cdot \IIm(\ch)$. Thus, the class \eqref{dRIntEta} is uniquely determined by $\hat{u}$ and $\hat{u}'$, up to the quotient by the subgroup $\Td(u) \cdot \IIm(\ch)$.
\begin{Def}\label{DiffClass} Let $\hat{u}$ and $\hat{u}'$ two differential Thom classes of $\pi: E \rightarrow X$ such that $I(\hat{u}) = I(\hat{u}')$ and $\Td(\hat{u}') = \Td(\hat{u})$. We define the \emph{difference class}:
	\[\delta(\hat{u}, \hat{u}') := \biggl[ \dR \biggl( \int_{E/X} \eta \biggr) \biggr] \in \frac{H^{-1}_{\dR}(X; \mathfrak{h}^{\bullet}_{\R})}{\Td(u) \cdot \IIm\, \ch}.
\]
\end{Def}
\begin{Theorem}\label{HomotDelta} With the hypotheses of definition \ref{DiffClass}, $\hat{u}$ and $\hat{u}'$ are homotopic if and only if $\delta(\hat{u}, \hat{u}') = 0$. Thus, the set of homotopy classes of differential Thom classes on a bundle $\pi: E \rightarrow X$, refining the same topological class $u$ and with a fixed Todd class $\Td(\hat{u})$, is a torsor over the group:
	\[\frac{H^{-1}_{\dR}(X; \mathfrak{h}^{\bullet}_{\R})}{\IIm\, \ch \cdot \Td(u)}.
\]
\end{Theorem}
For the proof see \cite[pp.\ 125-126]{Bunke}. The following corollary extends lemma \ref{Rule23Top} to the differential case.
\begin{Corollary}[2x3 principle]\label{Rule23Diff} Given two bundles $E, F \rightarrow X$, with projections $p_{E}: E \oplus F \rightarrow E$ and $p_{F}: E \oplus F \rightarrow F$, we consider a triple $(\hat{u}, \hat{v}, \hat{w})$ of differential Thom classes on $E$, $F$ and $E \oplus F$ respectively, such that $\hat{w}$ is homotopic to $p_{E}^{*}\hat{u} \cdot p_{F}^{*}\hat{v}$. Two elements of such a triple uniquely determine the third one up to homotopy.
\end{Corollary}
\begin{proof} Let us suppose that $\hat{w}$ and $\hat{u}$ are fixed up to homotopy. We choose a topological Thom class $v$ on $F$ such that $p_{E}^{*}u \oplus p_{E}^{*}v = w$. Let $\hat{v}'$ be any differential refinement of $v$ with Todd class $\Td(\hat{v}) = \Td(\hat{u})^{-1} \cdot \Td(\hat{w})$. The class $\hat{v}$ we are looking for must be of the form $\hat{v} := \hat{v}' - a(\eta)$, with $d\int_{E/X} \eta = 0$. The homotopy class of $\hat{v}$ is completely determined by $[\dR(\int_{E/X}\eta)]$. We have that $\delta(\hat{u}\hat{v}, \hat{w}) = \delta(\hat{u}\hat{v}', \hat{w}) - [\dR(\int_{E \oplus F/X}p_{F}^{*}R(\hat{u}) \cdot p_{E}^{*}\eta)] = \delta(\hat{u}\hat{v}', \hat{w}) - [\Td(u)][\dR(\int_{E/X}\eta)]$, hence the only possibility is $[\dR(\int_{E/X}\eta)] = [\Td(u)]^{-1} \cdot \delta(\hat{u}\hat{v}', \hat{w})$. \end{proof}

Finally, we will need the following corollary of theorem \ref{HomotDelta}.
\begin{Corollary}\label{CanOrientTriv} On the trivial bundle $X \times \R^{N}$ there is a canonical homotopy class of differential Thom classes, refining the canonical topological one.
\end{Corollary}
\begin{proof} Let us consider the line bundle $\{pt\} \times \R$. We put the canonical Thom class $u_{0}$, defined as the suspension of $1 \in \h^{0}$. For any refinement, the Todd class is necessary $1$. We choose $\hat{u}_{0}$ in such a way that $\int_{S^{1}} \hat{u}_{0} = 1$, being $S^{1} \simeq \R^{+}$. We show that $\hat{u}_{0}$ exists and is unique up to homotopy. Let $\hat{u}_{0}'$ be any differential refinement of $u_{0}$. Then $\int_{S^{1}} \hat{u}_{0}' = 1 + a(x)$, therefore we choose $\hat{u}_{0} := \hat{u}_{0}' - a(x dt)$. If $\hat{u}_{0}$ and $\hat{u}_{0}'$ both satisfy $\int_{S^{1}} \hat{u}_{0} = \int_{S^{1}} \hat{u}_{0}'$, then $\hat{u}_{0}' = \hat{u}_{0} + a(\eta)$ with $a(\int_{S^{1}} \eta) = 0$, i.e., $\int_{S^{1}} \eta \in \IIm(\ch)$. Hence, $\delta(\hat{u}_{0}, \hat{u}_{0}') = 0$, thus $\hat{u}_{0}$ is unique up to homotopy. On the bundle $\{pt\} \times \R^{N}$ we put the class $\hat{u}_{0, N} = \pi_{1}^{*}\hat{u}_{0} \cdots \pi_{N}^{*}\hat{u}_{0}$ and, on the bundle $X \times \R^{N}$, the class $\pi_{\R^{N}}^{*}\hat{u}_{0, N}$. \end{proof}

\subsection{Differential integration} The Gysin map $f_{!}: \hat{h}^{\bullet}(Y) \rightarrow \hat{h}^{\bullet - n}(X)$, being $n = \dim Y - \dim X$, is defined similarly to \eqref{GysinEmbeddingTop}, starting from the same data and refining the Thom class of the normal bundle to a differential one:
\begin{equation}\label{GysinEmbedding}
	f_{!}(\hat{\alpha}) = \int_{\R^{N}}i_{*}\varphi_{*}(\hat{u} \cdot \pi^{*}\hat{\alpha}).
\end{equation}
The integration map $\int_{\R^{N}}: \hat{h}^{\bullet+N}_{\cpt}(X \times \R^{N}) \rightarrow \hat{h}^{\bullet}(X)$ is defined as follows. There is a natural projection $\pi^{N}: (S^{1})^{N} \rightarrow S^{N}$, defined thinking of $S^{N}$ as $S^{1} \wedge \ldots \wedge S^{1} = (S^{1})^{N} / (S^{1} \vee \ldots \vee S^{1})$. For $\hat{\alpha} \in \hat{h}^{n+N}(X \times S^{N})$, we define:
\begin{equation}\label{IntSn}
	\int_{S^{n}}\hat{\alpha} := \int_{S^{1}} \cdots \int_{S^{1}} (\id \times \pi^{N})^{*}\hat{\alpha}.
\end{equation}
Given a class $\hat{\alpha} \in \hat{h}^{n+N}_{\cpt}(X \times \R^{N})$, since $S^{N}$ is the one-point compactification of $\R^{N}$, considering the embedding $j: \R^{N} \hookrightarrow S^{N}$ we can naturally define $j_{*}\hat{\alpha} \in \hat{h}^{n+N}(X \times S^{N})$ and $\int_{\R^{N}} \hat{\alpha} := \int_{S^{N}} j_{*}\hat{\alpha}$. We remark that, for any $N_{1}, N_{2}$ such that $N_{1} + N_{2} = N$, we have $\int_{\R^{N}} \hat{\alpha} = \int_{\R^{N_{2}}} \int_{\R^{N_{1}}} \hat{\alpha}$. That's because, being $\hat{\alpha}$ compactly supported, its extensions commute with any map between compactifications of $\R^{N}$; since $(S^{1})^{N} \simeq (S^{1})^{N_{2}} \times (S^{1})^{N_{1}}$, we get the result.

We define a \emph{representative of an $\hat{h}^{\bullet}$-orientation of $f$} as in definition \ref{DiffOrientedMap}, but considering a differential Thom class. Fixing such a representative $(\iota, \hat{u}, \varphi)$, the Gysin map $f_{!}$ is well-defined. Moreover, there is a natural map on differential forms, called \emph{curvature map}:
	\[\begin{split}
	R_{(\iota, \hat{u}, \varphi)}: \; & \Omega^{\bullet}(Y; \h^{\bullet}_{\R}) \rightarrow \Omega^{\bullet-n}(X; \h^{\bullet}_{\R}) \\
	& \omega \mapsto \int_{X \times \R^{N}/X} i_{*}\varphi_{*}(R(\hat{u}) \wedge \pi^{*}\omega).
\end{split}\]
It is easy to prove from the axioms that:
\begin{equation}\label{CurvAMaps}
	R(f_{!}\hat{\alpha}) = R_{(\iota, \hat{u}, \varphi)}(R(\hat{\alpha})) \qquad f_{!}a(\omega) = a(R_{(\iota, \hat{u}, \varphi)}(\omega)).
\end{equation}
The following definition is analogous to \ref{HomotopyOrientations}, but it keeps into account the curvature map. We define \emph{proper} representatives as in the topological framework.
\begin{Def}\label{HomotopyDiffOrientations} A \emph{homotopy} between two representatives $(\iota, \hat{u}, \varphi)$ and $(\iota', \hat{u}', \varphi')$ of an $\hat{h}^{\bullet}$-orientation of $f: Y \rightarrow X$ is a representative $(J, \hat{U}, \Phi)$ of an $\hat{h}^{\bullet}$-orientation of $\id \times f: I \times Y \rightarrow I \times X$, such that:
\begin{itemize}
	\item $(J, I(\hat{U}), \Phi)$ is proper over a neighborhood $V \subset I$ of $\{0, 1\}$;
	\item $(J, \hat{U}, \Phi)\vert_{f_{0}} = (\iota, \hat{u}, \varphi)$ e $(J, \hat{U}, \Phi)\vert_{f_{1}} = (\iota', \hat{u}', \varphi')$;
	\item $\pi_{X}^{*} \circ R_{(\iota, \hat{u}, \varphi)} = R_{(J, \hat{U}, \Phi)} \circ \pi_{Y}^{*}$.
\end{itemize}
\end{Def}
In particular, it follows that $R_{(\iota, \hat{u}, \varphi)} = R_{(\iota', \hat{u}', \varphi')}$. The following lemma is easy to prove.
\begin{Lemma} Let $(\iota, \hat{u}, \varphi)$ e $(\iota, \hat{u}', \varphi)$ be two \emph{proper} representatives such that $\hat{u}$ and $\hat{u}'$ are homotopic as differential Thom classes. Then the two representatives are homotopic.
\end{Lemma}
Thanks to corollary \ref{CanOrientTriv}, we define the equivalence of representatives up to stabilization as in the topological framework (def.\ \ref{EquivStab}).
\begin{Def} A \emph{$\hat{h}^{\bullet}$-orientation} on $f: Y \rightarrow X$ is an equivalence class $[\iota, \hat{u}, \varphi]$ of representatives, up to the equivalence relation generated by homotopy and stabilization.
\end{Def}
As a consequence of formula \eqref{HomFormula}, $f_{!}$ only depends on the $\hat{h}^{\bullet}$-orientation of $f$, not on the specific representative \cite[sec.\ 4.10]{Bunke}. We now consider a submersion $f: Y \rightarrow X$. In this case the Gysin map provides a good notion of integration.
\begin{Lemma}\label{SubPropDiff} Let $f: Y \rightarrow X$ be a submersion. Then, for any neat embedding $\iota: Y \hookrightarrow X \times \R^{N}$ and any Thom class $\hat{u}$ of the normal bundle, there exists a \emph{proper} representative $(\iota, \hat{u}, \varphi)$ of an $\hat{h}^{\bullet}$-orientation of $f$. Moreover, any two such representatives are homotopic, hence we can denote the orientation by $[\iota, \hat{u}]$.
\end{Lemma}
\begin{proof} The proof of the first statement is the same as theorem \ref{SubProp}. Because of the uniqueness up to homotopy of the tubular neighborhood and the embedding, we can find a homotopy $(\id \times \iota, \pi^{*}_{I}\hat{u}, \Phi)$ between the two representatives, which is proper. We have:
	\[\begin{split}
	R_{(\id \times \iota, \pi^{*}_{I}\hat{u}, \Phi)}(\pi_{Y}^{*}\omega) &= \int_{I \times X \times \R^{N}/I \times X} I_{*}\Phi_{*}(\pi^{*}_{I}R(\hat{u}) \cdot \pi^{*}\omega) \\
	&= \pi_{X}^{*} \int_{N_{\iota(Y)}(X \times \R^{N})} R(\hat{u}) \cdot \pi^{*}\omega = \pi_{X}^{*}(R_{(\iota, \hat{u}, \varphi)}(\omega)).
\end{split}\]
\end{proof}

\begin{Lemma}\label{ProperMorfMod} If the representative $(\iota, \hat{u}, \varphi)$ is proper, then the Gysin map is a morphism of $\hat{h}^{\bullet}(X)$-modules, i.e., for every $\hat{\alpha} \in \hat{h}^{\bullet}(Y)$ and $\hat{\beta} \in \hat{h}^{\bullet}(X)$:
	\[f_{!}(\hat{\alpha} \cdot f^{*}(\hat{\beta})) = f_{!}(\hat{\alpha}) \cdot \hat{\beta}.
\]
\end{Lemma}
\begin{proof} 
\begin{align*}
	f_{!}(\hat{\alpha} \cdot f^{*}\hat{\beta}) & = \int_{\R^{N}}\iota_{!}(\hat{\alpha} \cdot f^{*}\hat{\beta}) = \int_{\R^{N}} i_{*}\varphi_{*} (\hat{u} \cdot \pi^{*}\hat{\alpha} \cdot \pi^{*} \iota^{*} \pi_{X}^{*} \hat{\beta}) \\
	& = \int_{\R^{N}} i_{*}\varphi_{*} (\hat{u} \cdot \pi^{*}\hat{\alpha} \cdot \pi^{*} (\pi_{X}\vert_{\iota(Y)})^{*} \hat{\beta}) \\
	&= \int_{\R^{N}} i_{*} \bigl( \varphi_{*}(\hat{u} \cdot \pi^{*}\hat{\alpha}) \cdot (\varphi_{*}\pi^{*} (\pi_{X}\vert_{\iota(Y)})^{*} \hat{\beta}) \bigr) \\
	& \overset{\eqref{PropRepres}}= \int_{\R^{N}} i_{*} \bigl( \varphi_{*}(\hat{u} \cdot \pi^{*}\hat{\alpha}) \cdot (\pi_{X}\vert_{U}^{*} \hat{\beta}) \bigr) \\
	& = \int_{\R^{N}} i_{*} \bigl( \varphi_{*}(\hat{u} \cdot \pi^{*}\hat{\alpha}) \bigr) \cdot \pi_{X}^{*}\hat{\beta} = \biggl( \int_{\R^{N}}\iota_{!}(\hat{\alpha}) \biggr) \cdot \hat{\beta} = f_{!}(\hat{\alpha}) \cdot \hat{\beta}.
\end{align*}
\end{proof}
\begin{Lemma}\label{IntProp} If $(\iota, \hat{u}, \varphi)$ is proper, then:
	\[R_{(\iota, \hat{u}, \varphi)}(\omega) = \int_{Y/X} \Td(\hat{u}) \wedge \omega.
\]
\end{Lemma}
\begin{proof} Since the representative is proper, the tubular neighborhood $U$ is fibered over $\iota(Y)$, with fibers contained in $\R^{N}$; moreover, the disjoint union of the fibers over $\iota(f^{-1}\{x\})$ is a fiber over $x$. Hence, integrating with respect to $X \times \R^{N}/X$ a form, whose support is contained in $U$, is equivalent to integrating over $U/\iota(Y)$ and then over $Y/X$. Hence, we have:
	\[\begin{split}
	R_{(\iota, \hat{u}, \varphi)}\omega &= \int_{X \times \R^{N}/X} i_{*}\varphi_{*}(R(\hat{u}) \wedge \pi^{*}\omega) = \int_{Y/X} \int_{U/\iota(Y)} \varphi_{*}(R(\hat{u}) \wedge \pi^{*}\omega) \\
	&= \int_{Y/X} \biggl( \int_{N_{\iota(Y)}/\iota(Y)} R(\hat{u}) \biggr) \wedge \omega = \int_{Y/X} \Td(\hat{u}) \wedge \omega.
\end{split}\]
\end{proof}
\begin{Corollary}\label{ProperRA} If $(\iota, \hat{u}, \varphi)$ is proper, then:
	\[R(f_{!}\hat{\alpha}) = \int_{Y/X} \Td(\hat{u}) \wedge R(\hat{\alpha}) \qquad\qquad f_{!}a(\omega) = a \biggl( \int_{Y/X} \Td(\hat{u}) \wedge \omega \biggr).
\]
\end{Corollary}
\begin{proof} It immediately follows from lemma \ref{IntProp}. \end{proof}

\begin{Lemma}\label{DiffOrientedMapComposition} Let $f: Y \rightarrow X$ and $g: X \rightarrow W$ be $h^{\bullet}$-oriented neat submersions. There is a naturally induced $h^{\bullet}$-orientation on $g \circ f: Y \rightarrow W$ such that $(g \circ f)_{!} = g_{!} \circ f_{!}$.
\end{Lemma}
\begin{proof} We call $[\iota, u]$ and $[\kappa, v]$ the orientations of $f$ and $g$, and we construct an orientation $[\xi, w]$ on $g \circ f$ as in definition \ref{OrientedMapComposition}. The reader can prove that it does not depend on the representatives chosen for $[\iota, u]$ and $[\kappa, v]$. It remains to show that $(g \circ f)_{!} = g_{!} \circ f_{!}$. In order to prove this, we choose tubular neighborhoods $U$ and $V$ and diffeomorphisms $\varphi: N_{\iota(Y)}(X \times \R^{N}) \rightarrow U$ and $\psi: N_{\kappa(X)}(W \times \R^{L}) \rightarrow V$, in such a way that the representatives $(\iota, u, \varphi)$ and $(\kappa, v, \psi)$ are proper. We define the proper representative $(\xi, w, \nu)$ in the following way:
\begin{itemize}
	\item in order to define the tubular neighborhood of $\xi(Y)$ in $W \times \R^{L+N}$, for each point of the image of $\varphi$ in $X \times \R^{N}$, we consider the image under $(\psi, \id_{\R^{N}})$ of the corresponding fiber of $N_{\kappa(X) \times \R^{N}}(W \times \R^{L+N}) \simeq \pi^{*}_{N}N_{\kappa(X)}(W \times \R^{L}) \simeq N_{\kappa(X)}(W \times \R^{L}) \times \R^{N}$;
	\item $\nu$ is defined in the following way: given a vector $(A, B)_{y} \in N_{\iota(Y)}(X \times \R^{N}) \oplus (\pi^{*}_{N}N_{\kappa(X)}(W \times \R^{L}))\vert_{\xi(Y)}$, we apply $\varphi$ to $A$ in $y$ getting a point $p \in \{f(y)\} \times \R^{N}$; since $\pi^{*}_{N}N_{\kappa(X)}(W \times \R^{L})$, restricted to $\{f(y)\} \times \R^{N}$, coincides with $(N_{\kappa(X)}(W \times \R^{L}))_{\kappa(f(y))} \times \R^{N}$, we apply $(\psi, 1_{\R^{N}})$ to $B$, the latter translated from $\xi(y)$ to $(f(y), p)$, and we get $\nu(A, B)$.
\end{itemize}
Basically the parallel translation along $\R^{N}$ allows us to identify the bundle $(N_{\kappa(X) \times \R^{N}}(W \times \R^{L+N}))\vert_{\kappa(U)}$, being $U$ the image of $\varphi$, with the ``propagation'' along $U$ of $(N_{\kappa(X) \times \R^{N}}(W \times \R^{N+L}))\vert_{\xi(Y)}$. By definition, the orientation of $(N_{\kappa(X) \times \R^{N}}(W \times \R^{L+N}))\vert_{\kappa(U)}$ coincides with the ``propagation'' of the orientation of $(N_{\kappa(X) \times \R^{N}}(W \times \R^{N+L}))\vert_{\xi(Y)}$. Let us consider the following diagram:
\begin{equation}\label{DiagramGysinMaps}
	\xymatrix{
	\hat{h}^{\bullet}(Y) \ar[rr]^(.37){\iota_{!}} \ar[drr]_{f_{!}} \ar@/^2pc/[rrrr]^{((\kappa,\id_{\R^{N}}) \circ \iota)_{!}} & & \hat{h}^{\bullet+N-n}_{\cpt}(X \times \R^{N}) \ar[rr]^(.45){(\kappa,\id_{\R^{N}})_{!}} \ar[d]^{\int_{\R^{N}}} & & \hat{h}^{\bullet+N+L-n-l}_{\cpt}(W \times \R^{L+N}) \ar[d]^{\int_{\R^{N}}} \ar@/^5pc/[dd]^{\int_{\R^{L+N}}} \\
	& & \hat{h}^{\bullet-n}(X) \ar[rr]^(.45){\kappa_{!}} \ar[drr]_{g_{!}} & & \hat{h}^{\bullet+L-n-l}_{\cpt}(W \times \R^{L}) \ar[d]^{\int_{\R^{L}}} \\
	& & & & \hat{h}^{\bullet-n-l}(W).
}
\end{equation}
The curved arrows define $(g \circ f)_{!}$, therefore we have to prove the commutativity of the diagram. The left and lower triangles commute by definition. Commutativity of the right cell is a property of integration that we already remarked. About the upper cell, one has:\footnote{In the third line of the following equation we should have $\pi^{*}_{N}\hat{v}$, without the restriction to $V$, since $V$ does not appear in the definition of $(\kappa,\id_{\R^{N}})_{!}$. Hence, we are actually calculating $((\kappa,\id_{\R^{N}})\vert_{V})_{!}$. Nevertheless, since $\iota_{!}\hat{\alpha}$ vanishes outside $V$, the restriction is immaterial.}
	\[\begin{split}
	((\kappa,\id_{\R^{N}}) \circ \iota)_{!}\hat{\alpha} &= \psi_{\nu}^{*}\nu_{*}(\hat{\alpha} \cdot p_{1}^{*}\hat{u} \cdot p_{2}^{*}((\pi_{N}^{*}\hat{v})\vert_{\xi(Y)})) \\
	& = \psi_{\pi^{*}_{N}\psi}^{*}(\pi_{N}^{*}\psi)_{*}\bigl(\psi_{\varphi}^{*}\varphi_{*}(\pi^{*}\hat{\alpha} \cdot \hat{u}) \cdot (\pi_{N}^{*}\hat{v})\vert_{\xi(Y)})\bigr) \\
	& = \psi_{\pi^{*}_{N}\psi}^{*}(\pi_{N}^{*}\psi)_{*}(\pi^{*}\iota_{!}\alpha \cdot (\pi^{*}_{N}\hat{v})\vert_{V}) \\
	& = (\kappa,\id_{\R^{N}})_{!} \iota_{!} \alpha.
\end{split}\]
Finally, about the central square:
	\[\begin{split}
	\int_{\R^{N}} (\kappa,\id_{\R^{N}})_{!} \hat{\alpha} &= \int_{\R^{N}} (\pi_{N}^{*}j)_{*}(\pi_{N}^{*}\psi)_{*}(\pi^{*}\hat{\alpha} \cdot \pi_{N}^{*}\hat{v}) = j_{*}\psi_{*}\biggl(\int_{\R^{N}} \pi^{*}\hat{\alpha} \cdot \pi_{N}^{*}\hat{v}\biggr) \\
	&= j_{*}\psi_{*}\biggl(\biggl(\int_{\R^{N}} \pi^{*}\hat{\alpha} \biggr) \cdot \pi_{N}^{*}\hat{v}\biggr) = \kappa_{!}\biggl(\int_{\R^{N}} \hat{\alpha} \biggr).
\end{split}\]
\end{proof}
The following lemma is a consequence of lemma \ref{Rule23Diff} and of the uniqueness up to homotopy and stabilization of the embedding $\iota$.
\begin{Lemma}[2x3 principle]\label{Rule23Subm} Let $f: Y \rightarrow X$ and $g: X \rightarrow W$ be $\hat{h}^{\bullet}$-oriented neat submersions, with orientations $[\iota, \hat{u}]$ and $[\kappa, \hat{v}]$, and let $[\xi, \hat{w}]$ be the orientation induced on $g \circ f$, as shown in lemma \ref{DiffOrientedMapComposition}. Two elements of the triple $([\iota, \hat{u}], [\kappa, \hat{v}], [\xi, \hat{w}])$ uniquely determines the third one.
\end{Lemma}
Finally, we consider differential orientations of manifolds. We define them as in the topological case (def.\ \ref{OrientedManifold} and \ref{OrientedManifoldBoundary}). We put $\Td(X) := \Td(\hat{u})$, being $\hat{u}$ the orientation of the stable normal bundle. If follows from corollary \ref{ProperRA} that:
	\[R((p_{X})_{!}(\hat{\alpha})) = \int_{X} \Td(X) \wedge R(\hat{\alpha}).
\]
Let us consider a submersion $f: Y \rightarrow X$ between compact $\hat{h}^{\bullet}$-oriented manifolds. Since $p_{Y} = p_{X} \circ f$, it follows from lemma \ref{Rule23Subm} that $f$ inherits a unique orientation from the ones of $X$ and $Y$. Hence, the integration map $f_{!}: \hat{h}^{\bullet}(Y) \rightarrow \hat{h}^{\bullet-n}(X)$, for $n = \dim \, Y - \dim \, X$, is well-defined for submersions between compact $\hat{h}^{\bullet}$-oriented manifolds.
\begin{Lemma}\label{ARoofBoundary} For $X$ an $\hat{h}^{\bullet}$-manifold with boundary and $\Phi: X \rightarrow I$ defined as in \ref{OrientedManifoldBoundary}, we have:
	\[\int_{0}^{1}R(\Phi_{!}\hat{\alpha}) = \int_{X} \Td(X) \wedge R(\hat{\alpha}).
\]
\end{Lemma}
\begin{proof} Because of formula \eqref{CurvAMaps}, we have:
	\[\begin{split}
	\int_{0}^{1} R(\Phi_{!}\hat{\alpha}) &= \int_{0}^{1} \int_{I \times \R^{N}/I} i_{*}\varphi_{*}(R(\hat{u}) \wedge \pi^{*}R(\hat{\alpha})) = \int_{N_{\iota(X)}(I \times \R^{N})} R(\hat{u}) \wedge \pi^{*}R(\hat{\alpha}) \\
	&= \int_{X} \biggl( \int_{N_{\iota(X)}(I \times \R^{N})/\iota(X)} R(\hat{u}) \biggr) \wedge R(\hat{\alpha}) = \int_{X} \Td(X) \wedge R(\hat{\alpha}).
\end{split}\]
\end{proof}
About the restriction of an orientation to the boundary, the same considerations we did in the topological framework hold for differential extensions. In particular, formula \eqref{RestrictionBoundaryGysin} keeps on holding.

\section{Cohomology and homology}\label{CohHom}

Let us consider a multiplicative cohomology theory $h^{\bullet}$ represented by a spectrum $(E_{\bullet}, e_{\bullet}, \varepsilon_{\bullet})$, where $e_{n}$ is the marked point of $E_{n}$ and $\varepsilon_{n}: (\Sigma E_{n}, \Sigma e_{n}) \rightarrow (E_{n+1}, e_{n+1})$ is the structure map, whose adjoint $\tilde{\varepsilon}_{n}: E_{n} \rightarrow \Omega_{e_{n+1}} E_{n+1}$ is a homotopy equivalence. Considering the spectrum $(E_{\bullet} \wedge X, e_{\bullet} \wedge x_{0}, \varepsilon_{\bullet} \wedge 1)$, the dual homology theory $h_{\bullet}$ is defined, on a space with marked point $(X, x_{0})$, as \cite{Whitehead}:
\begin{equation}\label{DualHomology}
	h_{n}(X, x_{0}) := \pi_{n}(E_{\bullet} \wedge X, e_{\bullet} \wedge x_{0}) = \varinjlim_{k} \pi_{n+k}(E_{k} \wedge X, e_{k} \wedge x_{0}).
\end{equation}
The unreduced groups are defined as $h_{n}(X) := h_{n}(X_{+}, \infty)$, for $X_{+} = X \sqcup \{\infty\}$. There is a natural map for every $n \in \Z$ \cite[pp.\ 289-290]{Switzer}:
\begin{equation}\label{HomCohomSpectra}
	\xi^{n}: h^{n}(X) \rightarrow \Hom_{\h^{\bullet}}(h_{n-\bullet}(X), \h^{\bullet}).
\end{equation}
From \eqref{HomCohomSpectra} we can easily define:
\begin{equation}\label{HomCohomRSpectra}
	\xi^{n}_{\mathbb{R}}: h^{n}(X) \otimes_{\mathbb{Z}} \mathbb{R} \rightarrow \Hom_{\mathfrak{h}^{\bullet}}(h_{n-\bullet}(X), \mathfrak{h}^{\bullet}_{\mathbb{R}}).
\end{equation}
It follows from the universal coefficient theorem (see \cite[prop.\ 13.5 p.\ 285]{Adams} and \cite[p.\ 290]{Switzer}) that the map \eqref{HomCohomRSpectra} is an isomorphism. Finally, for $\mathfrak{h}^{\bullet}_{\mathbb{R}/\mathbb{Z}} := h^{\bullet}(\{*\}; \mathbb{R}/\mathbb{Z})$, there is a natural map:
\begin{equation}\label{HomCohomRZSpectra}
	\xi^{n}_{\mathbb{R}/\mathbb{Z}}: h^{n}(X; \mathbb{R}/\mathbb{Z}) \rightarrow \Hom_{\mathfrak{h}^{\bullet}}(h_{n-\bullet}(X), \mathfrak{h}^{\bullet}_{\mathbb{R}/\mathbb{Z}}).
\end{equation}
For singular cohomology or K-theory \eqref{HomCohomRZSpectra} is an isomorphism, as a consequence of the universal coefficient theorem \cite{Yosimura, Lott}. This depends on the fact that ordinary cohomology and K-theory are Pontrjagin self-dual \cite{FMS}.

In \cite{Jakob} the author provides a geometric construction of the homology theory dual to a given cohomology theory, which we briefly recall in the following, only in the case of a single space $X$. With respect to \cite{Jakob}, in the following definition we replace the quotient by diffeomorphisms and vector bundle modification with the quotient by the Gysin map associated to a submersion. Actually, we could not require that the map is a submersion, but this hypothesis will make it easier to define the differential extension of cycles.
\begin{Def} Let $h^{\bullet}$ be a multiplicative cohomology theory. On a space $X$ with the homotopy type of a finite CW-complex, we define:
\begin{itemize}
	\item the group of \emph{$n$-precycles of $h_{\bullet}$} as the free abelian group generated by the quadruples $(M, u, \alpha, f)$, with:
\begin{itemize}
	\item $(M, u)$ a smooth compact manifold (without boundary) with $h^{\bullet}$-orientation $u$, whose connected components $\{M_{i}\}$ have dimension $n+q_{i}$, with $q_{i}$ arbitrary;
	\item $\alpha \in h^{\bullet}(M)$, such that $\alpha\vert_{M_{i}} \in h^{q_{i}}(M)$;
	\item $f: M \rightarrow X$ a continuous map;
\end{itemize}
	\item the group of \emph{$n$-cycles of $h_{\bullet}$}, denoted by $z_{n}(X)$, as the quotient of the group of $n$-precycles by the free subgroup generated by elements of the form:
\begin{itemize}
	\item $(M, u, \alpha + \beta, f) - (M, u, \alpha, f) - (M, u, \beta, f)$;
	\item $(M, u, \alpha, f) - (M_{1}, u\vert_{M_{1}}, \alpha\vert_{M_{1}}, f\vert_{M_{1}}) - (M_{2}, u\vert_{M_{2}}, \alpha\vert_{M_{2}}, f\vert_{M_{2}})$, for $M = M_{1} \sqcup M_{2}$;
	\item $(M, u, \varphi_{!}\alpha, f) - (N, v, \alpha, f \circ \varphi)$ for $\varphi: N \rightarrow M$ a submersion and $\varphi_{!}: h^{\bullet}(N) \rightarrow h^{\bullet}(M)$ the Gysin map;
\end{itemize}
	\item the group of \emph{$n$-boundaries of $h_{\bullet}$}, denoted by $b_{n}(X)$, as the subgroup of $z_{n}(X)$ generated by the cycles which are representable by a pre-cycle $(M, u, \alpha, f)$, such that there exits a quadruple $(W, U, A, F)$, where $W$ is a manifold and $M = \partial W$, $U$ is an $h^{\bullet}$-orientation of $W$ and $U\vert_{M} = u$, $A \in h^{\bullet}(W)$ and $A\vert_{M} = \alpha$, $F: W \rightarrow X$ is a continuous map satisfying $F\vert_{M} = f$.
\end{itemize}
We define $h_{n}(X) := z_{n}(X) / b_{n}(X)$.
\end{Def}
Describing in this way the dual homology theory, the map \eqref{HomCohomSpectra} corresponds to:
\begin{equation}\label{HomCohom}
\begin{split}
	\xi^{n}: \; & h^{n}(X) \rightarrow \Hom_{\h^{\bullet}}(h_{n-\bullet}(X), \h^{\bullet}) \\
	& \alpha \mapsto \bigl( [M, u, \beta, f] \mapsto (p_{M})_{!}(\beta \cdot f^{*}\alpha) \bigr),
\end{split}
\end{equation}
where $p_{M}: M \rightarrow \{pt\}$. We verify that \eqref{HomCohom} is well-defined. If we consider a function between compact manifolds $\varphi: N \rightarrow M$ and two representatives $(M, u, \varphi_{!}\beta, f)$ and $(N, v, \beta, f \circ \varphi)$ of the homology class, we have:
	\[\begin{split}
	\xi^{n}(\alpha)[N, v, \beta, f \circ \varphi] = (p_{N})_{!}&(\beta \cdot \varphi^{*}f^{*}\alpha) = (p_{M})_{!} \varphi_{!}(\beta \cdot \varphi^{*}f^{*}\alpha) \\
	&= (p_{M})_{!} (\varphi_{!}\beta \cdot f^{*}\alpha) = \xi^{n}(\alpha)[M, u, \varphi_{!}\beta, f].
\end{split}\]
Let us now suppose that $(M, u, \beta, f) = \partial (W, U, B, F)$. Then we consider a function $\Phi: W \rightarrow [0, 1]$, such that $\Phi^{-1}(0) = M$ and $\Phi^{-1}(1) = \emptyset$. We have that $(p_{M})_{!}(\beta \cdot f^{*}\alpha) = \Phi_{!}(B \cdot F^{*}\alpha)\vert_{\{0\}}$. By homotopy invariance, such a class coincides with $\Phi_{!}(B \cdot F^{*}\alpha)\vert_{\{1\}} = 0$. Finally, the image of $\alpha$ is a $\h^{\bullet}$-module homomorphism, since, for $\gamma \in \h^{t}$:
	\[\begin{split}
	\xi^{n}(\alpha)([(M, u, \beta, f)] \cdot \gamma) = \xi^{n}(\alpha)[M, u, \beta \cdot &(p_{M})^{*}\gamma, f] = (p_{M})_{!}(\beta \cdot f^{*}\alpha \cdot (p_{M})^{*}\gamma) \\
	& = (p_{M})_{!}(\beta \cdot f^{*}\alpha) \cdot \gamma = \xi^{n}(\alpha)[M, u, \beta, f] \cdot \gamma.
\end{split}\]
Tensorizing with $\R$, we get the isomorphism:
\begin{equation}\label{HomCohomR}
	\xi^{n}_{\R}: h^{n}(X) \otimes_{\Z} \R \overset{\simeq}\longrightarrow \Hom_{\h^{\bullet}}(h_{n-\bullet}(X), \h^{\bullet}_{\R}).
\end{equation}
The map \eqref{HomCohomRZSpectra} corresponds to:
\begin{equation}\label{HomCohomRZ1}
\begin{split}
	\xi^{n}_{\R/\Z}: \; & h^{n}(X; \R/\Z) \rightarrow \Hom_{\h^{\bullet}}(h_{n-\bullet}(X), \h^{\bullet}_{\R/\Z}) \\
	& \alpha \mapsto \bigl( [M, u, \beta, f] \mapsto (p_{M})_{!}(\beta \cdot f^{*}\alpha) \bigr).
\end{split}
\end{equation}
The product $\beta \cdot f^{*}\alpha$ is provided by the structure of $h^{\bullet}$-module on $h^{\bullet}(\,\cdot\,; \R/\Z)$.

\section{Flat pairing}\label{FlatPairing}

\subsection{Flat classes} Given a smooth map $f: Y \rightarrow X$, with $n = \dim Y - \dim X$, the Gysin map $f_{!}: \hat{h}^{\bullet}(Y) \rightarrow \hat{h}^{\bullet - n}(X)$ previously defined depends on the $\hat{h}^{\bullet}$-orientation, but, if we restrict to flat classes, it only depends on the topological $h^{\bullet}$-orientation. In order to prove this statement, we show that there is a natural graded module structure on $\hat{h}^{\bullet}_{\fl}(X)$ over $h^{\bullet}(X)$, i.e.\ there exists a product:
\begin{equation}\label{ModuleFlat}
	h^{\bullet}(X) \otimes_{\Z} \hat{h}^{\bullet}_{\fl}(X) \rightarrow \hat{h}^{\bullet}_{\fl}(X).
\end{equation}
This is due to the fact that the product of differential classes restricts to a product $\hat{h}^{\bullet}(X) \otimes_{\Z} \hat{h}^{\bullet}_{\fl}(X) \rightarrow \hat{h}^{\bullet}_{\fl}(X)$, since, being the curvature multiplicative, if one of the two factors has vanishing curvature, also the result has. Moreover, the product $\hat{\alpha} \cdot \hat{\beta}$, with $\hat{\beta}$ flat, only depends on $I(\hat{\alpha})$. In fact, if $I(\hat{\alpha}) = 0$, then $\hat{\alpha} = a(\omega)$. Because of definition \ref{MultDiffExt}, we have $a(\omega) \cdot \hat{\beta} = a(\omega \wedge R(\hat{\beta})) = 0$. We can show in the same way that also the product $\hat{h}^{\bullet}_{\cpt}(E) \otimes_{\Z} \hat{h}^{\bullet}(E) \rightarrow \hat{h}^{\bullet}_{\cpt}(E)$ can be refined to $h^{\bullet}_{\cpt}(E) \otimes_{\Z} \hat{h}^{\bullet}_{\fl}(E) \rightarrow \hat{h}^{\bullet}_{\fl,\cpt}(E)$. Therefore, given a real vector bundle $E \rightarrow X$ of rank $n$ with (topological) Thom class $u$, we define the Thom morphism:
	\[\begin{split}
	T_{\fl}: \; & \hat{h}^{\bullet}_{\fl}(X) \rightarrow \hat{h}^{\bullet + n}_{\fl,\cpt}(E) \\
	& \hat{\alpha} \mapsto u \cdot \pi^{*}\hat{\alpha}.
\end{split}\]
From this it easily follows that the Gysin map $f_{!}$, when applied to a flat class, only depends on the topological orientation. Lemmas \ref{GysinMapProp1} and \ref{GysinMapProp2} keep on holding, with the same proofs. The following lemma is analogous to \ref{ProperRA}, but it is not necessary to suppose that the orientation is proper.
\begin{Lemma}\label{GysinC10} For $f: Y \rightarrow X$ a map of $h^{\bullet}$-oriented manifolds and $\theta \in H^{\bullet-1}_{\dR}(Y; \h^{\bullet}_{\R})$, we have:
	\[f_{!}(a(\theta)) = a(f_{!}(\Td(u) \cdot \theta)).
\]
This is equivalent to the fact that, for any $\alpha \in h^{\bullet}(X) \otimes_{\Z} \R$:\footnote{In equation \eqref{PushForwardChernC} we are considering the Chern character as defined on $h^{\bullet}(X) \otimes_{\Z} \R$, in which case it is an isomorphism. If we consider it as defined on $h^{\bullet}(X)$, then $a(\ch \alpha) = 0$, and formula \eqref{PushForwardChernC} implies coherently that $f_{!}(a(\ch\alpha)) = 0$.}
\begin{equation}\label{PushForwardChernC}
	f_{!}(a(\ch\,\alpha)) = a(\ch(f_{!}\alpha)).
\end{equation}
\end{Lemma}
\begin{proof} Let us consider a differential Thom class $\hat{u}$ of $N_{\iota(Y)}(X \times \R^{N})$ refining the orientation $u$ induced by the ones of $X$ and $Y$. We have:
	\[\begin{split}
	f_{!}(a(\theta)) &= i_{*}\varphi_{*}(\hat{u} \cdot \pi^{*}a(\theta)) = i_{*}\varphi_{*}(a(\dR(R(\hat{u})) \cdot \pi^{*}\theta) = a(i_{*}\varphi_{*}(\ch\,u \cdot \pi^{*}\theta)) \\
	&= a(i_{*}\varphi_{*}(\ch^{(n)}u \cdot \pi^{*}(\Td(u) \cdot \theta)) = a(f_{!}(\Td(u) \cdot \theta)).
\end{split}\]
Formula \eqref{PushForwardChernC} follows from the Grothendieck-Riemann-Roch theorem.
\end{proof}
\begin{Corollary} The Gysin map associated to $f: Y \rightarrow X$ induces a morphism of exact sequences of $\h^{\bullet}$-modules:
	\[\xymatrix{
	\cdots \ar[r] & h^{\bullet}(Y) \ar[r] \ar[d] & h^{\bullet}(Y) \otimes_{\Z} \R \ar[r] \ar[d] & \hat{h}^{\bullet+1}_{\fl}(Y) \ar[r] \ar[d] & h^{\bullet+1}(Y) \ar[r] \ar[d] & \cdots \\
	\cdots \ar[r] & h^{\bullet}(X) \ar[r] & h^{\bullet}(X) \otimes_{\Z} \R \ar[r] & \hat{h}^{\bullet+1}_{\fl}(X) \ar[r] & h^{\bullet+1}(X) \ar[r] & \cdots
}\]
where the map $h^{\bullet}(X) \otimes_{\Z} \R \rightarrow \hat{h}^{\bullet+1}_{\fl}(X)$ is defined by $\alpha \mapsto a(\ch \alpha)$.
\end{Corollary}

\subsection{Flat pairing}

We can now define a natural $\hat{\h}^{\bullet}_{\fl}$-valued pairing on a manifold $X$ between $\hat{h}^{\bullet}_{\fl}$ and $h_{\bullet}$, that, in the case of singular differential cohomology, reduces to the holonomy of a flat Deligne cohomology class. When $\hat{h}^{\bullet}_{\fl} \simeq h^{\bullet}(\,\cdot\,; \R/\Z)$, such a pairing coincides with formula \eqref{HomCohomRZ1}.
\begin{Def} For $X$ a differential manifold, there is a natural pairing:
\begin{equation}\label{FlatPairing2}
\begin{split}
	\xi^{n}_{\fl}: \; & \hat{h}^{n}_{\fl}(X) \rightarrow \Hom_{\mathfrak{h}^{\bullet}} (h_{n-\bullet}(X), \hat{\h}^{\bullet}_{\fl}) \\
	& \hat{\alpha} \mapsto \bigl( [M, u, \beta, f] \mapsto (p_{M})_{!}(\beta \cdot f^{*}\hat{\alpha}) \bigr).
\end{split}
\end{equation}
We recall that $p_{M}$ is the unique map $p_{M}: M \rightarrow \{pt\}$ and the product $\beta \cdot f^{*}\hat{\alpha}$ is defined by \eqref{ModuleFlat}. The invariance by $\h^{\bullet}$ is defined by:
\begin{equation}\label{PairingInvariance}
	\xi^{n}_{\fl}(\hat{\alpha})([M, u, \beta, f] \cdot \gamma) = \xi^{n}_{\fl}(\hat{\alpha})([M, u, \beta, f]) \cdot \gamma.
\end{equation}
\end{Def}
\vspace{0.2cm}
In order to show that \eqref{FlatPairing2} is well-defined, i.e.\ that it does not depend on the representative $(M, u, \beta, f)$ of the homology class, and that formula \eqref{PairingInvariance} holds, we use use an argument similar to the one used about \eqref{HomCohom}.
\begin{Lemma}\label{MorphismComplexes} There is a morphism of complexes of $\h^{\bullet}$-modules (the second one not being exact in general):
\begin{scriptsize}
	\[\xymatrix{
	\cdots \ar[r]^(.3){r} & h^{n}(X) \otimes_{\Z} \R \ar[r]^{a} \ar[d]^{\xi^{n}_{\R}} & \hat{h}^{n+1}_{\fl}(X) \ar[r]^{I} \ar[d]^{\xi^{n+1}_{\fl}} & h^{n+1}(X) \ar[r]^(.75){r} \ar[d]^{\xi^{n+1}} & \cdots \\
	\cdots \ar[r]^(.3){r'} & \Hom_{\h^{\bullet}}(h_{n-\bullet}(X), \h^{\bullet}_{\R}) \ar[r]^(.47){a'} & \Hom_{\h^{\bullet}}(h_{n+1-\bullet}(X), \hat{\h}^{\bullet}_{\fl}) \ar[r]^{I'} & \Hom_{\h^{\bullet}}(h_{n+1-\bullet}(X), \h^{\bullet}) \ar[r]^(.75){r'} & \cdots
}\]
\end{scriptsize}
\end{Lemma}
\begin{proof} We only have to prove the commutativity of the square under the map $a$. It easily follows from the fact that, for $\alpha \in h^{\bullet}(X) \otimes_{\Z} \R$ and $\beta \in h^{\bullet}(X)$:
	\[a(\ch\alpha) \cdot \beta = a(\ch(\alpha\beta)).
\]
That's because, for any differential refinement $\hat{\beta}$ of $\beta$, we have $a(\ch\alpha) \cdot \hat{\beta} = a(\ch\alpha \cdot R(\hat{\beta})) = a(\ch\alpha \cdot \ch\beta) = a(\ch(\alpha\beta))$. \end{proof}

We call $\h^{n}_{\Z}$ the image of the Chern character $\ch: \h^{n} \rightarrow H^{n}_{\dR}(pt; \h^{\bullet}_{\R}) \simeq \h^{n}_{\R}$, which coindices with $\alpha \mapsto \alpha \otimes_{\Z} \R$.
\begin{Theorem}\label{NoTorsionRZ} If $\h^{\bullet}$ has no torsion, the pairing \eqref{FlatPairing2} is an isomorphism and $\hat{\h}^{\bullet}_{\fl} \simeq \h^{\bullet-1}_{\R}/\h^{\bullet-1}_{\Z}$.
\end{Theorem}
\begin{proof} The isomorphism $\hat{\h}^{\bullet}_{\fl} \simeq \h^{\bullet-1}_{\R}/\h^{\bullet-1}_{\Z}$ easily follows from the long exact sequence $\cdots \rightarrow \h^{\bullet-1} \rightarrow \h^{\bullet-1}_{\R} \rightarrow \h^{\bullet-1}_{\fl} \rightarrow \h^{\bullet} \rightarrow \cdots$. In order to prove that \eqref{FlatPairing2} is an isomorphism, we consider the functor:
	\[k^{n}(X) := \Hom_{\mathfrak{h}^{\bullet}} (h_{n-\bullet}(X), \h^{\bullet-1}_{\R}/\h^{\bullet-1}_{\Z}),
\]
defined on the same category on which $h^{\bullet}$ is defined. The functor $\hat{h}^{\bullet}_{\fl}$ is a cohomology theory \cite[Sec.\ 7]{BS}. We show that also $k^{\bullet}$ is a cohomology theory. In fact, homotopy invariance follows from the one of $h_{\bullet}$. Moreover, since $\h^{\bullet-1}_{\R}/\h^{\bullet-1}_{\Z}$ is a division graded abelian group, the functor $\Hom_{\mathfrak{h}^{\bullet}} (\,\cdot\,, \h^{\bullet-1}_{\R}/\h^{\bullet-1}_{\Z})$ is exact, hence $k^{\bullet}$ has a natural Mayer-Vietoris sequence \cite[sec.\ 7]{BS}. Finally, we show that $\xi^{\bullet}_{\fl}(pt): \hat{\h}^{\bullet}_{\fl} \rightarrow k^{\bullet}(pt)$ is an isomorphism. Let us consider $1 \in \h^{0}$ and the class $u_{1} := [pt, u_{0}, 1, \id] \in h_{0}(pt)$. Any class belonging to $h_{n}(pt)$ is of the form $[pt, u_{0}, \gamma, \id] = u_{1} \cdot \gamma$. By $\h^{\bullet}$-invariance, for any $\hat{\alpha} \in \h^{n}_{\fl}$, we have $\xi^{n}_{\fl}(\hat{\alpha})(u_{1} \cdot \gamma) = \xi^{n}_{\fl}(\hat{\alpha})(u_{1}) \cdot \gamma$, hence the morphism $\xi^{n}_{\fl}(\hat{\alpha})$ is completely determined by the image of $u_{1}$, belonging to $\hat{\h}^{n}_{\fl}$. It follows that $\xi^{n}_{\fl}: \hat{\h}^{n}_{\fl} \rightarrow k^{n}(pt)$ is an isomorphism. \end{proof}

\section{Generalized Cheeger-Simons characters}\label{DiffCycles}

We now describe a model of the homology groups $h_{\bullet}$, that involves a differential extension of the cycles, so that we will be able to define the generalized Cheeger-Simons characters. 
\begin{Def} On a smooth compact manifold $X$, we define:
\begin{itemize}
	\item the group of \emph{$n$-precycles} as the free abelian group generated by the quadruples $(M, \hat{u}, \hat{\alpha}, f)$, with:
\begin{itemize}
	\item $(M, \hat{u})$ a smooth compact manifold without boundary with $\hat{h}^{\bullet}$-orientation $\hat{u}$, whose connected components $\{M_{i}\}$ have dimension $n+q_{i}$, with $q_{i}$ arbitrary;
	\item $\hat{\alpha} \in \hat{h}^{\bullet}(M)$, such that $\hat{\alpha}\vert_{M_{i}} \in \hat{h}^{q_{i}}(M)$;
	\item $f: M \rightarrow X$ a smooth map;
\end{itemize}
	\item the group of \emph{$n$-cycles}, denoted by $\hat{z}_{n}(X)$, as the quotient of the group of $n$-precycles by the free subgroup generated by elements of the form:
\begin{itemize}
	\item $(M, \hat{u}, \hat{\alpha} + \hat{\beta}, f) - (M, \hat{u}, \hat{\alpha}, f) - (M, \hat{u}, \hat{\beta}, f)$;
	\item $(M, \hat{u}, \hat{\alpha}, f) - (M_{1}, \hat{u}\vert_{M_{1}}, \hat{\alpha}\vert_{M_{1}}, f\vert_{M_{1}}) - (M_{2}, \hat{u}\vert_{M_{2}}, \hat{\alpha}\vert_{M_{2}}, f\vert_{M_{2}})$, for $M = M_{1} \sqcup M_{2}$;
	\item $(M, \hat{u}, \varphi_{!}\hat{\alpha}, f) - (N, \hat{v}, \hat{\alpha}, f \circ \varphi)$ for $\varphi: N \rightarrow M$ a submersion, oriented via the 2x3 principle;
\end{itemize}
	\item the group of \emph{$n$-boundaries}, denoted by $\hat{b}_{n}(X)$, as the subgroup of $\hat{z}_{n}(X)$ generated by the cycles which are representable by a pre-cycle $(M, \hat{u}, \hat{\alpha}, f)$ such that there exists a quadruple $(W, \hat{U}, \hat{A}, F)$, where $W$ is a manifold and $M = \partial W$, $\hat{U}$ is an $\hat{h}^{\bullet}$-orientation of $W$ and $\hat{U}\vert_{M} = \hat{u}$, $\hat{A} \in \hat{h}^{\bullet}(W)$ such that $\hat{A}\vert_{M} = \hat{\alpha}$, and $F: W \rightarrow X$ is a smooth map satisfying $F\vert_{M} = f$.
\end{itemize}
We define $h'_{n}(X) := \hat{z}_{n}(X) / \hat{b}_{n}(X)$.
\end{Def}
\begin{Theorem}\label{TopDiffIso} The natural group morphism:
	\[\begin{split}
	\Phi: \, &h'_{\bullet}(X) \rightarrow h_{\bullet}(X) \\
	& [(M, \hat{u}, \hat{\alpha}, f)] \rightarrow [(M, I(\hat{u}), I(\hat{\alpha}), f)]
\end{split}\]
is an isomorphism.
\end{Theorem}
\begin{proof} We divide the proof in three steps.
\subparagraph{\emph{Step 1.}} If $I(\hat{u}) = I(\hat{u}')$ and $I(\hat{\alpha}) = I(\hat{\alpha}')$, then $[(M, \hat{u}, \hat{\alpha}, f)] = [(M, \hat{u}', \hat{\alpha}', f)]$ in $h'_{\bullet}(X)$. In fact, since $\hat{\alpha}' = \hat{\alpha} + a(\rho)$, we consider on $I \times M$ the class $A = \pi_{I}^{*}\hat{\alpha} + a(t\cdot \pi_{I}^{*}\rho)$, which links $\alpha$ to $\hat{\alpha}'$. Moreover, we orient the projection $I \times M \rightarrow I$ in the following way: the embedding $M \hookrightarrow \R^{N}$ naturally determines an embedding $I \times M \hookrightarrow I \times \R^{N}$, with normal bundle $\pi_{I}^{*}(N_{M}\R^{N})$; we put on such a bundle the Thom class $\hat{U} = \pi_{I}^{*}\hat{u} + a(t\cdot \pi_{I}^{*}\eta)$. Hence $\partial(I \times M, \hat{U}, \hat{A}, \id \times f) = (M, \hat{u}, \hat{\alpha}, f) - (M, \hat{u}', \hat{\alpha}', f)$.
\subparagraph{\emph{Step 2.}} Given two equivalent topological precycles $(M, u, \varphi_{!}\alpha, f) \simeq (N, v, \alpha, f \circ \varphi)$, any two differential refinements $(M, \hat{u}, \hat{\alpha}', f)$ (with $I(\hat{\alpha}') = \varphi_{!}\alpha$) and $(N, \hat{v}, \hat{\alpha}, f \circ \varphi)$ are equivalent in $h'_{\bullet}(X)$. In fact, by definition $[(M, \hat{u}, \varphi_{!}\hat{\alpha}, f)] = [(N, \hat{v}, \hat{\alpha}, f \circ \varphi)]$. By the first step, this implies that $[(M, \hat{u}, \hat{\alpha}', f)] = [(N, \hat{v}, \hat{\alpha}, f \circ \varphi)]$.
\subparagraph{\emph{Step 3.}} The  morphism $\Phi$ is clearly well-defined and surjective. Therefore, we only have to prove injectivity. Let us suppose that $\Phi[(M, \hat{u}, \hat{\alpha}, f)] = 0$. Then $[(M, I(\hat{u}), I(\hat{\alpha}), f)]$ is equivalent, as a cocycle, to $[(N, v, \beta, g)]$ such that $(N, v, \beta, g) = \partial (W, V, B, G)$. This means that there exists a sequence of pre-cycles $(M_{i}, u_{i}, \alpha_{i}, f_{i})$, for $i = 0, \ldots, n$, such that $(M_{0}, u_{0}, \alpha_{0}, f_{0}) = (M, I(\hat{u}), I(\hat{\alpha}), f)$, $(M_{n}, u_{n}, \alpha_{n}, f_{n}) = (N, v, \beta, g)$ and such that there exists a submersion $\varphi_{i}: M_{i} \rightarrow M_{i+1}$ or $\psi_{i}: M_{i+1} \rightarrow M_{i}$ such that $f_{i} = f_{i+1} \circ \varphi_{i}$ and $\alpha_{i+1} = \varphi_{!}(\alpha_{i})$, or the analogue for $\psi_{i}$. We choose a differential refinement $(M_{i}, \hat{u}_{i}, \hat{\alpha}_{i}, f_{i})$ for each $i$, such that for $i = 0$ it coincides with $(M, \hat{u}, \hat{\alpha}, f)$. By the second step, we get that $[(M_{i}, \hat{u}_{i}, \hat{\alpha}_{i}, f_{i})] = [(M_{i+1}, \hat{u}_{i+1}, \hat{\alpha}_{i+1}, f_{i+1})]$, hence $[(M, \hat{u}, \hat{\alpha}, f)] = [(N, \hat{v}, \hat{\beta}, g)]$. We now consider a differential refinement $(W, \hat{V}, \hat{B}, G)$ of $(W, V, B, G)$. By the step 1, $[(N, \hat{v}, \hat{\beta}, g)] = [(N, \hat{V}\vert_{N}, \hat{B}\vert_{N}, g)] = 0$. \end{proof}

\begin{Def}\label{GeneralizedCS} A \emph{Cheeger-Simons differential $\hat{h}^{\bullet}$-character} of degree $n$ on $X$ is a couple $(\chi_{n}, \omega_{n})$, where:
\begin{equation}\label{GeneralizedCSDef}
	\chi_{n} \in \Hom_{\hat{\h}^{\bullet}} (\hat{z}_{n-\bullet}(X), \hat{\h}^{\bullet}) \quad\qquad \omega_{n} \in \Omega^{n}(X; \h^{\bullet}_{\R})
\end{equation}
such that, if $(M, \hat{u}, \hat{\beta}, f) = \partial(W, \hat{U}, \hat{B}, F)$, then:
\begin{equation}\label{CSFormula}
	\chi_{n}[M, \hat{u}, \hat{\beta}, f] = -a \biggl( \int_{W} \Td(W) \wedge R(\hat{B}) \wedge F^{*}\omega_{n} \biggr).
\end{equation}
The $\hat{\h}^{\bullet}$-invariance is defined by:
\begin{equation}\label{CSInvariance}
	\chi_{n}(\hat{\alpha})([M, \hat{u}, \hat{\beta}, f] \cdot \gamma) = \chi_{n}(\hat{\alpha})[M, \hat{u}, \hat{\beta}, f] \cdot \gamma.
\end{equation}
We denote by $\check{h}^{n}(X)$ the group of characters of degree $n$.
\end{Def}

We briefly comment formula \eqref{CSFormula}. Let us suppose that $[M, \hat{u}, \hat{\beta}, f] \in \hat{z}_{n-k}(X)$ and that $M$ is connected. Then $\dim(M) = n - k + q$ and $\hat{\beta} \in \hat{h}^{q}(M)$, hence $\dim(W) = n - k + q + 1$ and $\hat{B} \in \hat{h}^{q}(M)$. Thus, in the r.h.s.\ of \eqref{CSFormula}, we integrate on $W$ a $\h^{\bullet}_{\R}$-valued form of degree $n + q + 0$, hence we get a form on the point of degree $n + q - (n - k + q + 1) = k - 1$. Applying $a$, we get a class belonging to $\hat{\h}^{k}$, as desired.
\begin{Theorem}\label{CShnThm} There is a natural graded-group morphism:
\begin{equation}\label{CShn}
\begin{split}
	CS_{\hat{h}}^{\bullet}: \; & \hat{h}^{\bullet}(X) \rightarrow \check{h}^{\bullet}(X) \\
	& \hat{\alpha} \mapsto (\chi, R(\hat{\alpha})),
\end{split}
\end{equation}
where $\chi$ is defined, for $[M, \hat{u}, \hat{\beta}, f] \in \hat{z}_{n-k}(X)$, by:
	\[\chi[M, \hat{u}, \hat{\beta}, f] := (p_{M})_{!}(\hat{\beta} \cdot f^{*}\hat{\alpha}).
\]
\end{Theorem}
\begin{proof} If we consider two representatives $(M, u, \varphi_{!}\beta, f)$ and $(N, v, \beta, f \circ \varphi)$ of the homology class, we have, thanks to lemmas \ref{ProperMorfMod} and \ref{DiffOrientedMapComposition}:
	\[\begin{split}
	\chi[N, \hat{v}, \hat{\beta}, f \circ \varphi] = (p_{N})_{!}&(\hat{\beta} \cdot \varphi^{*}
	f^{*}\hat{\alpha}) = (p_{M})_{!} \varphi_{!}(\hat{\beta} \cdot \varphi^{*}f^{*}\hat{\alpha}) \\
	&= (p_{M})_{!} (\varphi_{!}\hat{\beta} \cdot f^{*}\hat{\alpha}) = \chi[M, \hat{u}, \varphi_{!}\hat{\beta}, f].
\end{split}\]
Let us now suppose that $(M, \hat{u}, \hat{\beta}, f) = \partial (W, \hat{U}, \hat{B}, F)$. Then, for $\Phi$ defined as in \ref{OrientedManifoldBoundary}, thanks to formula \eqref{RestrictionBoundaryGysin} one has:
	\[(p_{M})_{!}(\hat{\beta} \cdot f^{*}\hat{\alpha}) = (\Phi_{!}(\hat{B} \cdot F^{*}\hat{\alpha}))\vert_{\{0\}}.
\]
Since $(\Phi_{!}(\hat{B} \cdot F^{*}\hat{\alpha}))\vert_{\{1\}} = 0$, because $\Phi^{-1}(1) = \emptyset$, from the homotopy formula \eqref{HomFormula} and theorem \ref{ARoofBoundary} we have:
	\[(p_{M})_{!}(\hat{\beta} \cdot f^{*}\hat{\alpha}) = -a\biggl( \int_{I} R(\Phi_{!}(\hat{B} \cdot F^{*}\hat{\alpha})) \biggr) = -a\biggl( \int_{W} \Td(W) \wedge R(\hat{B} \cdot F^{*}\hat{\alpha}) \biggr).
\]
Hence:
	\[\chi[M, \hat{u}, \hat{\beta}, f] = -a\biggl( \int_{W} \Td(W) \wedge R(\hat{B}) \wedge F^{*}R(\hat{\alpha}) \biggr).
\]
This is exactly formula \eqref{CSFormula} for $\omega_{n} = R(\hat{\alpha})$. Finally:
	\[\begin{split}
	\chi(\hat{\alpha})([M, \hat{u}, \hat{\beta}, f] \cdot \hat{\gamma}) &= (p_{M})_{!}(f^{*}\hat{\alpha} \cdot \hat{\beta} \cdot (p_{M})^{*}\hat{\gamma}) \\
	&= (p_{M})_{!}(f^{*}\hat{\alpha} \cdot \hat{\beta}) \cdot \hat{\gamma} = \chi(\hat{\alpha})[M, \hat{u}, \hat{\beta}, f] \cdot \hat{\gamma}.
\end{split}\]
\end{proof}

The proof of the following theorem is straightforward from the previous definition.
\begin{Theorem}\label{HolonomyFlat} When $\hat{\alpha}$ is flat, the value of the associated Cheeger-Simons character over $[M, \hat{u}, \hat{\beta}, f]$ coincides with the value of \eqref{FlatPairing2} on the corresponding homology class.
\end{Theorem}
Considering the pairing \eqref{FlatPairing2}, we remark that there is an embedding:
	\[j: \Hom_{\mathfrak{h}^{\bullet}} (h_{n-\bullet}(X), \hat{\h}^{\bullet}_{\fl}) \hookrightarrow \check{h}^{n}(X).
\]
In fact, a morphism $\varphi_{n} \in \Hom_{\mathfrak{h}^{\bullet}} (h_{n-\bullet}(X), \hat{\h}^{\bullet}_{\fl})$ determines a unique morphism $\chi_{n}: \hat{z}_{n-\bullet}(X) \rightarrow \hat{\h}^{\bullet}$ defined by $\chi_{n}[M, \hat{u}, \hat{\beta}, f] := \varphi_{n}[M, I(\hat{u}), I(\hat{\beta}), f]$, and we define $j(\varphi_{n}) := (\chi_{n}, 0)$. It follows from formula \eqref{CSFormula} that the image of $j$ is the subgroup of generalized Cheeger-Simons characters with vanishing curvature, which we call $\check{h}^{n}_{\fl}(X)$. Let us consider the embedding $i: \hat{h}^{\bullet}_{\fl}(X) \hookrightarrow \hat{h}^{\bullet}(X)$. The following diagram commutes:
	\[\xymatrix{
	\hat{h}^{n}_{\fl}(X) \ar[r]^(.3){\xi^{n}_{\fl}} \ar[d]_{i} & \Hom_{\mathfrak{h}^{\bullet}} (h_{n-\bullet}(X), \hat{\h}^{\bullet}_{\fl}) \ar[d]^{j} \\
	\hat{h}^{n}(X) \ar[r]^{CS^{n}_{\hat{h}}} & \check{h}^{n}(X).
}\]
Therefore $i$ restricts to an embedding $i': \Ker(\xi_{\fl}^{n}) \hookrightarrow \Ker(CS^{n}_{\hat{h}})$, and $j$ restricts to an embedding $j': \IIm(\xi^{n}_{\fl}) \hookrightarrow \IIm(CS^{n}_{\hat{h}})$. Because of $j$ and $j'$ we can construct a morphism $a: \Coker(\xi_{\fl}^{n}) \rightarrow \Coker(CS^{n}_{\hat{h}})$. We now show that actually $i'$ and $a$ are isomorphisms. In particular, when $\h^{\bullet}$ has no torsion, $CS^{n}_{\hat{h}}$ is an isomorphism, because of theorem \ref{NoTorsionRZ}.
\begin{Theorem} The following canonical isomorphisms hold:
\begin{equation}\label{KerCokerEqual}
	\Ker(\xi_{\fl}^{n}) \simeq \Ker(CS^{n}_{\hat{h}}), \qquad \Coker(\xi_{\fl}^{n}) \simeq \Coker(CS^{n}_{\hat{h}}).
\end{equation}
\end{Theorem}
\begin{proof} If $\hat{\alpha} \in \hat{h}^{n}(X)$ is not flat, then $CS^{n}_{\hat{h}}(\hat{\alpha}) \neq 0$, since $CS^{n}_{\hat{h}}(\hat{\alpha}) = (\chi_{n}, R(\hat{\alpha}))$ and $R(\hat{\alpha}) \neq 0$. Hence $\Ker(CS^{n}_{\hat{h}}) \subset \Ker(\xi_{\fl}^{n})$ and the equality follows. Moreover, $\check{h}^{n}_{\fl}(X) \cap \IIm(CS^{n}_{\hat{h}}) = \IIm(\xi_{\fl}^{n})$, hence $a: \Coker(\xi_{\fl}^{n}) \rightarrow \Coker(CS^{n}_{\hat{h}})$ is an embedding. If $(\chi_{n}, \omega_{n}) \in \check{h}^{n}(X)$, we consider a class $\hat{\alpha} \in \hat{h}^{n}(X)$ such that $R(\hat{\alpha}) = \omega_{n}$, and we call $(\chi'_{n}, \omega_{n}) := CS_{\hat{h}}(\hat{\alpha})$. Then $(\chi'_{n} - \chi_{n}, 0) \in \check{h}^{n}_{\fl}(X)$, and, in $\Coker(CS^{n}_{\hat{h}})$, one has $[(\chi_{n}, \omega_{n})] = [(\chi'_{n} - \chi_{n}, 0)] \in \IIm\,a$. Therefore $a$ is also surjective. \end{proof}

\subsection{Singular cohomology and K-theory}

In the case of singular cohomology, we have $\hat{\h}^{0} \simeq \Z$, $\hat{\h}^{1} \simeq \R/\Z$ and $\hat{\h}^{n} = 0$ for $n \neq 0,1$. Moreover, $\hat{\h}_{\fl}^{1} \simeq \R/\Z$ and $\hat{\h}_{\fl}^{n} = 0$ for $n \neq 1$. Thus, a flat differential class $\hat{\alpha} \in \hat{H}^{n}_{\fl}(X)$ defines a differential character $\chi_{n}: H_{n-1}(X) \rightarrow \mathbb{R}/\mathbb{Z}$. Such a character corresponds (up to the exponential) to the holonomy of the corresponding smooth Deligne cohomology class, for which there is an explicit formula \cite{GT2}. Actually, the holonomy is defined on singular cycles, with no need of differential refinement. Thus, there is a canonical isomorphism between the group of Cheeger-Simons characters defined on singular cycles (in the usual sense), and the group of Cheeger-Simons characters defined on differential cycles (as defined in the present paper), being both groups canonically isomorphic to $\hat{H}^{n}(X)$. We can explicitly describe this isomorphism for a large class of cycles. In fact, let us consider an $(n-1)$-cycle $[(M, \hat{u}, 1, f)$, for $\hat{u}$ any differential orientation refining the topological one. Then, for $\alpha \in \hat{H}^{n}(X)$, $f^{*}\hat{\alpha}$ is flat and with trivial first Chern class for dimensional reasons, therefore $f^{*}\hat{\alpha} = [h] \in H^{n}(M; \R)/H^{n}(M; \Z)$ and the holonomy coincides with the exponential of $h([M])$, being $[M]$ the fundamental class of $M$.

Finally, we consider the case of complex K-theory. We call $\kk^{\bullet}$ the $K$-theory ring of the point. In this case $I$ provides a canonical isomorphism:
\begin{equation}\label{IsoKK}
	\hat{\kk}^{2n} \simeq \kk^{2n} \quad \forall n \in \Z,
\end{equation}
because $\Omega^{2n-1}(pt, \kk^{\bullet}_{\R}) \simeq \kk^{2n-1}_{\R} = 0$. Moreover, let us consider the class $\gamma_{0} \in \kk^{-2}$, corresponding to the dual of the tautological line bundle of $\mathbb{P}^{1}(\C)$. Then Bott periodicity $K^{n}(X) \simeq K^{n-2}(X)$ is given by $\alpha \mapsto \alpha \cdot \gamma_{0}$. Such a periodicity can be extended to differential K-theory. In fact, because of \eqref{IsoKK}, there is a unique differential extension $\hat{\gamma}_{0}$, hence the map $\hat{\alpha} \mapsto \hat{\alpha} \cdot \hat{\gamma}_{0}$ is an isomorphism. It follows that $\hat{\kk}^{2n} \simeq \hat{\kk}^{0} \simeq \Z$. Moreover, $\hat{\kk}^{2n+1} \simeq \hat{\kk}^{1} \simeq \R/\Z$, because the elements of $\hat{\kk}^{1}$ are of the form $a(\theta)$, with $\theta \in H^{0}_{\dR}(pt; \kk^{\bullet}_{\R}) \simeq \kk^{0}_{\R} \simeq \R$, being $a(\theta) = 0$ if and only if $\theta \in \Z$. The periodicity can be extended to differential cycles, via the isomorphism:
	\[\begin{split}
	B: \;& \hat{z}_{n}(X) \rightarrow \hat{z}_{n-2}(X) \\
	&[M, \hat{u}, \hat{\beta}, f] \mapsto [M, \hat{u}, \hat{\beta} \cdot \hat{\gamma}_{0}, f].
\end{split}\]
Therefore, a generalized Cheeger-Simons character, as defined by formula \eqref{GeneralizedCSDef}, is uniquely determined by its restriction to $\hat{z}_{n-1}$, because of formula \eqref{CSInvariance}. The same holds for the pairing \eqref{FlatPairing2} (see \cite{Lott} for an analytic description of the pairing). It follows from theorem \ref{NoTorsionRZ} or from the universal coefficient theorem for K-theory \cite[formula 3.1]{Yosimura}, that \eqref{FlatPairing2} and \eqref{CShn} are isomorphisms. Therefore, in the case of K-theory, the pairing \eqref{FlatPairing2} and theorem \ref{CShnThm} can be summarized and enriched as follows.
\begin{Def} A \emph{Cheeger-Simons differential $\hat{K}^{\bullet}$-character} of degree $n$ on $X$ is a couple $(\chi_{n}, \omega_{n})$, where:
\begin{equation}
	\chi_{n}: \hat{z}_{n-1}(X) \rightarrow \mathbb{R}/\mathbb{Z} \qquad\quad \omega_{n} \in \Omega^{n}(X; \kk^{\bullet}_{\mathbb{R}})
\end{equation}
such that, if $(M, \hat{u}, \hat{\beta}, f) = \partial(W, \hat{U}, \hat{B}, F)$, then $\chi_{n}[(M, \hat{u}, \hat{\beta}, f)] \equiv_{\Z} \int_{W} \Td(W) \wedge R(\hat{B}) \wedge F^{*}\omega_{n}$. We denote by $\check{K}^{n}(X)$ the group of characters of degree $n$.
\end{Def}
\begin{Theorem} There is a natural group \emph{isomorphism}:
\begin{equation}
\begin{split}
	CS_{K}^{n}: \; & \hat{K}^{n}(X) \rightarrow \check{K}^{n}(X) \\
	& \hat{\alpha} \rightarrow (\chi, R(\hat{\alpha})),
\end{split}
\end{equation}
where $\chi$ is defined, for $[(M, \hat{u}, \hat{\beta}, f)] \in \hat{z}_{n-1}(X)$, by $\chi[(M, \hat{u}, \hat{\beta}, f)] := (p_{M})_{!}(\hat{\beta} \cdot f^{*}\hat{\alpha})$. Restricting to flat classes, we get an \emph{isomorphism}:
\begin{equation}
	\xi_{\fl}^{n}: \hat{K}^{n}_{\fl}(X) \rightarrow \Hom(K_{n-1}(X), \mathbb{R}/\mathbb{Z}).
\end{equation}
\end{Theorem}
This construction is equivalent to the one considered in \cite{BM}. In fact, let us call $\mathcal{C}_{\bullet}(X)$ the groups of $K$-cycles, as defined in \cite{BM}. An element of $\mathcal{C}_{\bullet}(X)$ is represented by a triple $(M, (E, h, \nabla), \phi)$, where $M$ is a compact manifold with fixed metric and spin$^{c}$-structure, $E \rightarrow X$ is a vector bundle with Hermitian metric $h$ and compatible connection $\nabla$, and $\phi: M \rightarrow X$ is a smooth function.\footnote{In \cite{BM} the function $\phi$ is only supposed to be continuous, but the result is the same by smooth approximation theorem.} Using the Freed-Lott model \cite[def.\ 2.15 p.\ 8]{FL}, $(E, h, \nabla)$ defines a differential $K$-theory class $\hat{\alpha}$. Moreover, the fixed metric and spin$^{c}$-structure define a Thom class of the tangent bundle of $M$ \cite[theorem C.12 p.\ 388]{LM}, hence, by the 2x3 principle, define an orientation of the stable normal bundle of $M$. In this way we get a differential $K$-cycle as defined in the present paper. Since the vector bundle modification, considered in \cite{BM}, is a particular case of the Gysin map, we get a morphism $\rho: \mathcal{C}_{\bullet}(X) \rightarrow \hat{z}_{\bullet}(X)$. Finally, let us call $\check{\mathcal{K}}^{\bullet}(X)$ the group of differential $K$-characters, as defined in \cite{BM}. We define $\nu^{\bullet}: \check{K}^{\bullet}(X) \rightarrow \check{\mathcal{K}}^{\bullet}(X)$ as $\nu(\chi)(A) := \chi(\rho(A))$ and the identity on the curvature. Both the groups $\check{K}^{\bullet}(X)$ and $\check{\mathcal{K}}^{\bullet}(X)$ fit in the exact sequence $0 \rightarrow K^{\bullet}(X; \R/\Z) \rightarrow \check{K}^{\bullet}(X) \rightarrow \Omega^{\bullet+1}_{0}(X) \rightarrow 0$ \cite[Theorem 1 p.\ 431]{BM}, and the morphism $\nu^{\bullet}$ commutes with the identities of $K^{\bullet}(X; \R/\Z)$ and $\Omega^{\bullet+1}_{0}(X)$. Because of the five lemma, $\nu^{\bullet}$ is an isomorphism.


\end{document}